\documentclass{article}
\usepackage[%
journal=JFG,    
lang=american,   
]{ems-journal}

\usepackage[all]{xy}

\theoremstyle{plain}
\newtheorem{theorem}{Theorem}[section]
\newtheorem{lemma}[theorem]{Lemma}
\newtheorem{proposition}[theorem]{Proposition}
\newtheorem{corollary}[theorem]{Corollary}
\theoremstyle{definition}
\newtheorem{example}[theorem]{Example}
\newtheorem{remark}[theorem]{Remark}

\newcommand{\bs}[1]{\mbox{$\boldsymbol{#1}$}}
\newcommand{\ind}{\mathrm{index}}

\numberwithin{equation}{section}

\begin{document}

\title{A combinatorial Fredholm module on self-similar sets built on $n$-cubes}
\titlemark{Fredholm module on \texorpdfstring{$n$}{n}-cube}

\emsauthor{1}{Takashi Maruyama}{T.~Maruyama}
\emsauthor{2}{Tatsuki Seto}{T.~Seto}


\emsaffil{1}{Takashi Maruyama, System Platform Research Laboratories, NEC Corporation, 
1753 Shimonumabe, Nakahara-ku, Kawasaki, Kanagawa, Japan \email{49takashi@nec.com} \email{49takashi@gmail.com}}
\emsaffil{2}{General Education and Research Center, Meiji Pharmaceutical University, 
2-522-1 Noshio, Kiyose-shi, Tokyo, Japan \email{tatsukis@my-pharm.ac.jp}}

\classification[28A80]{46L87}


\begin{abstract}
We construct a Fredholm module 
on self-similar sets such as 
the Cantor dust, the Sierpinski carpet 
and the Menger sponge. 
Our construction is 
a higher dimensional analogue of 
Connes' combinatorial construction of 
the Fredholm module on the Cantor set. 
We also calculate the Dixmier trace of two 
operators induced by the Fredholm module. 
\end{abstract}

\maketitle


\section*{Introduction}

In the 1990s, 
A. Connes \cite[Chapter IV]{MR1303779} 
introduced the quantized calculus based on the Fredholm modules. 
A Fredholm module on an involutive algebra $\mathcal{A}$ 
is a pair $(H, F)$ of a Hilbert space $H$ and bounded operator $F$ 
such that $\mathcal{A}$ acts on $H$ and
$a(F - F^{\ast}), a(F^{2}-1), [F,a] \in \mathcal{K}(H)$ 
for any $a \in \mathcal{A}$. 
The commutator $[F,a]$ is called a quantized differential of $a$. 
The notion and calculus of Fredholm modules provide many techniques in studying    various spaces. Such examples are noncompact spaces, foliated spaces, noncommutative spaces, and fractal spaces, to name a few. In the present paper, we study Fredholm modules on a special class of fractal spaces called self-similar sets. 

The first study of quantized calculus on self-similar sets 
is given by Connes \cite[Chapter IV]{MR1303779}. 
Connes defined the Fredholm module $(H,F)$ on $C(CS)$,  
where $CS$ is the Cantor set 
realized in the interval $[0,1]$, 
by using vertices of the removed intervals. 
Specifically, he set $H_{I} = \ell^{2}(\{ a \}) \oplus \ell^{2}(\{ b \})$ for an open interval $I = (a,b)$ and $F_{I} = \begin{bmatrix} & 1 \\ 1 & \end{bmatrix}$ on $H_{I}$, and constructed $(H,F)$ by taking a direct sum 
of $(H_{I}, F_{I})$ on all removed open intervals which appear along the construction of $CS$. 
The Fredholm module $(H,F)$ defines an element in $K^{0}(C(CS))$. 
Connes also 
calculated the non-vanishing Dixmier trace 
$\mathrm{Tr}_{\omega}(|[F,x]|^{\dim_{H}(CS)})$. Here, $x$ is the coordinate function on $\mathbb{R}$ 
(we consider $x$ as a multiplication operator) and 
$\dim_{H}(CS)$ is the Hausdorff dimension of $CS$.  
We call $|[F,x]|^{\dim_{H}(CS)}$ the quantized volume measure on $CS$
and 
$\mathrm{Tr}_{\omega}(|[F,x]|^{\dim_{H}(CS)})$ the quantized volume on $CS$ since the commutator $[F,x]$ is a quantized differential of $x$. 

In this paper, we generalize the quantized volume measure 
and the quantized volume to higher dimensional self-similar sets. 
For the generalization, we construct a higher dimensional analogue of the Connes' Fredholm module. 
We now present what we mean by 
the generalization of 
Connes' quantized volume 
$\mathrm{Tr}_{\omega}(|[F,x]|^{\dim_{H}(CS)})$. 
When we have the Fredholm module $(H_{K},F_{K})$ 
on a fractal set $K \subset \mathbb{R}^{n}$ 
such that an algebra of functions on $\mathbb{R}^{n}$ acts on $H_{K}$, 
a commutator of operators $[F_{K} , x^{\alpha}]$ ($\alpha = 1, \dots , n$) for the $\alpha$-th coordinate function $x^{\alpha}$ on $\mathbb{R}^{n}$ is obtained. 
The commutator $[F_{K} , x^{\alpha}]$ is a quantized differential of $x^{\alpha}$, 
hence we say the operator 
\[
|[F_{K} , x^{1}][F_{K} , x^{2}] \cdots [F_{K} , x^{n}]|^{p}
\]
is a quantized volume measure of the volume measure $dx^{1} dx^{2} \cdots dx^{n}$ 
on $\mathbb{R}^{n}$. Here, $p \in \mathbb{R}$ is defined by a fractal dimension on $K$. Then the value 
\[
\mathrm{Tr}_{\omega}(|[F_{K} , x^{1}][F_{K} , x^{2}] \cdots [F_{K} , x^{n}]|^{p})
\] 
may be also called as a generalization of Connes' quantized volume on $K$. 

\begin{figure}[h]
  \includegraphics[width = 120mm]{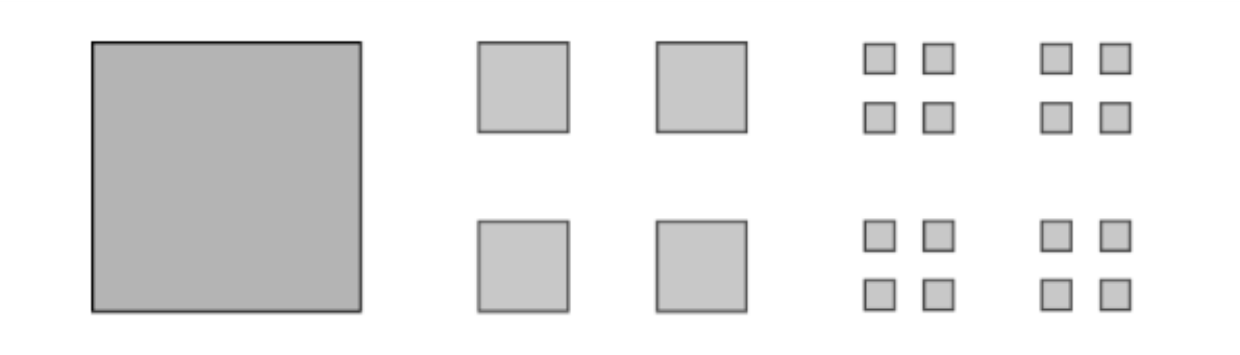}
\caption{The first $3$ steps of  construction of $\mathcal{CD}_{2}$. }
\label{fig:dust}
\end{figure}

Let us explain some examples that motivate us to conduct this work. 
We construct a Fredholm module on 
a self-similar set $K$ built on the square in $\mathbb{R}^{2}$. 
When we adopt a standard way to construct Fredholm modules on more generic $K$ (see \cite[section 2]{MR3743228}), 
it suffices to choose a subset $S \subset K$ for $(H_{S}, F_{S})$. 
As constructed in \cite{MR2238285}, 
when we choose $S = \{ a,b  \}$ ($2$ points), 
we have the same Fredholm module 
$(H_{S}, F_{S})$ as Connes' one. 
The Fredholm module $(H_{S}, F_{S})$ gives rise to a Fredholm module $(H_{K},F_{K})$ 
composed by the direct sum over all steps in the construction of $K$. 
Then the commutator $[F_{K}, x]$ (resp. $[F_{K}, y]$) is essentially given by 
the length of the projection of a segment $ab$ to the $x$-axis (resp. $y$-axis), and we can calculate the value 
$\mathrm{Tr}_{\omega}(|[F_{K} , x] [F_{K} , y]|^{p})$. However, the value may vanish: suppose that the vertices of the square are numbered counterclockwise in the order $v_{0} , v_{1} , v_{2}, v_{3}$. When $K$ is the Cantor dust (see Figure \ref{fig:dust}) and every edge of the square is parallel to either $x$- or $y$-axis,   
we have $[F_{K} , x] [F_{K} , y] = 0$ if 
$S = \{ v_{i} , v_{j} \}$ is the boundary of an edge of the square. 
On the other hand, if we choose 
$S = \{ v_{0}, v_{2} \}$ to be the boundary of a diagonal line of the square, 
we have the non-trivial value 
$\mathrm{Tr}_{\omega}(|[F_{K} , x][F_{K} , y]|^{\dim_{H}(K)/2})$.   
Therefore the subset $S = \{ v_{0}, v_{2} \}$ may look like an appropriate choice for the Cantor dust. However,  the value 
$\mathrm{Tr}_{\omega}(|[F_{K} , x][F_{K} , y]|^{\dim_{H}(K)/2})$ 
for $S = \{ v_{0}, v_{2} \}$ 
is not preserved under the rotation of the square. 
In fact, for a self-similar set $K$ obtained by rotation of the Cantor dust with rotation angle $\pi/4$ around $v_{0}$, 
we have $[F_{K} , x][F_{K} , y] = 0$ for $S = \{ v_{0}, v_{2} \}$.   
Thus, the choice of $S = \{a, b\}$ giving a non-trivial $\mathrm{Tr}_{\omega}(|[F_{K} , x][F_{K} , y]|^{p})$ depends on $K$. 
In this paper, we also present a way to construct a Fredholm module 
for $K$ that specifies a unified choice of $S$ (not necessarily $2$ points) 
and show that 
the Fredholm module induces a non-trivial higher dimensional quantized volume measure 
which is invariant under the Euclidean isometries in $\mathbb{R}^{n}$.

The outline of our construction of the Fredholm module $(H_{K} , F_{K})$ on $K$ is the following. 
Let $\gamma_{n} = [0,1]^{n}$ be the $n$-cube and 
$\{ f_{s} : \gamma_{n} \to \gamma_{n} \}$ ($s = 1,2,\dots , N$)
be similitudes with the similarity ratio 
$0 < r_{s} < 1$. 
We note that we do not require the open set condition. 
We now have a decreasing sequence of compact sets 
$\displaystyle K_{j} = \bigcup_{(s_{1}, \dots , s_{j})}f_{s_{1}} \circ \dots \circ f_{s_{j}}(\gamma_{n})$ in which each $f_{s_{1}} \circ \dots \circ f_{s_{j}}(\gamma_{n})$ 
is a small copy of the $n$-cube. Then, the sequence gives rise to the limiting set 
$\displaystyle K = \bigcap_{j=0}^{\infty}K_{j}$. 
Our construction of $(H_{K} , F_{K})$ is made of $2$ steps:  
the first step is the construction of 
the Fredholm module $(\mathcal{H},F_{n})$ on the $n$-cube; 
see subsection \ref{subsec:dfn-n}. 
In our construction, 
we use \textit{all} vertices (instead of $2$ points) of $n$-cubes, that is, 
we set $\mathcal{H} = \ell^{2}(\{ \text{vertices} \})$ 
with a suitable $\mathbb{Z}_{2}$-grading. 
In the definition of $F_{n}$, we use induction on the dimension $n$. 
The resulting Fredholm module represents 
the Kasparov product ($n$-times) of 
Connes' Fredholm module on an interval.  
The second step is taking the direct sum 
of $(\mathcal{H},F_{n})$ on 
all the copies of $n$-cubes; 
see subsection \ref{subsec:FMST}.  
Our Fredholm module 
$(\mathcal{H}_{K} , F_{K})$ 
is defined over $C(V_{K})$, 
where we denote by $V_{K}$ 
the closure of the vertices of 
all $n$-cubes $f_{s_{1}} \circ \dots \circ f_{s_{j}}(\gamma_{n})$. 
Note that, in general, $V_{K}$ includes $K$ properly, but $V_{K}$ coincides with $K$ for some important examples such as the Cantor dust, the Sierpinski carpet and the Menger sponge. 
Dividing by the length of edges of each $n$-cubes, we get the Dirac operator $D_{K}$ on $K$ and the spectral triple on $K$. 

Main results in the paper are basically twofold: our first result is the construction of a higher dimensional analogue of the Connes’ Fredholm module. 
This Fredholm module is also non-trivial in $K^{0}$ group under additional assumptions, which is given in Theorem \ref{thm:suffnonvanish} as a part of other properties of the Fredholm module delved in Section \ref{sec:FonK}. 
The second result is the derivation of concrete values for higher dimensional variants of the quantized volume measure and the quantized volume for some self-similar sets. The results are given in Section \ref{sec:Dixmier}. 
The calculation is based on a Clifford algebra's relation which the commutators $[F_{n} , x^{\alpha}]$ ($\alpha = 1, \dots , n$) generally satisfy for the $\alpha$-th coordinate functions $x^{\alpha}$ on $\mathbb{R}^{n}$. 
The Clifford algebra's relation is quantizaion of the relation of the exterior differentials $dx^{\alpha}$; see Proposition \ref{prp:qK} and \ref{prp:absqf} for the details.

Fredholm modules on self-similar sets are constructed by various researchers and studied from various aspects. F. Cipirani-J. Sauvageot \cite{MR2472035} constructed Fredholm modules on post critically finite fractals (p.c.f fractals) by regular harmonic structures. M. Ionescu-L. Rogers-A. Teplyaev \cite{MR2964679} studied weakly summable Fredholm modules  
in the cases of some finitely and infinitely ramified fractals. 
As an unbounded picture of Fredhlm modules, spectral triples on some self-similar sets have been also  extensively investigated.
E. Christensen-C. Ivan-L. Lapidus \cite{MR2357322} 
defined a spectral triple on the Sierpinski gasket 
$\mathcal{SG}$, which in turn defines an element in $K^{1}(C(\mathcal{SG}))$, by using the Dirac operator on the circle. D. Guido-T. Isola \cite{MR2238285} 
defined a spectral triple on self-similar sets with the open set condition in higher dimension 
by using Connes' Fredholm module on an interval. 
Guido-Isola \cite{MR3743228} also defined 
a spectral triple on nested fractals by using Connes' Fredholm module on an interval.   
See Introduction in \cite{MR3743228} for more related literatures.  

Let us compare our spectral triple with Guido-Isola's triples.  
First, our Fredholm module can not be constructed on self-similar sets on arbitrary subsets in $\mathbb{R}^{n}$, but on $n$-cubes. Our construction also does not require the open set condition. An example of the case for a self-similar set without the open set condition is given in subsection \ref{subsec:nonOSC}. Second, our triple and the triple in \cite{MR2238285} are not constructed on the algebra $C(K)$ of the continuous functions on $K$.  
Our algebra $C(V_{K})$ coincides with $C(K)$ for 
some important examples such as the Cantor dust, the Sierpinski carpet and the Menger sponge. The calculation of the value $\mathrm{Tr}_{\omega}(|D_{K}|^{-p})$ for our Dirac operator is also given in subsection \ref{subsec:Dixv}. 
The triple in \cite{MR3743228} is defined on $C(K)$ for the class of nested fractals, but the examples mentioned above are not the case. 
\begin{center} 
\begin{tabular}{c|c|c|c}
 & Ours & G-I's \cite{MR2238285} & G-I's \cite{MR3743228} \\\hline
space & self-similar set on $n$-cube & self-similar set on  $\mathbb{R}^{n}$ & nested fractal \\\hline 
algebra & $C(V_{K})$ & $C(C)$ & $C(K)$ \\ 
\end{tabular}
\end{center}

We will study 
more noncommutative geometry of our Fredholm module $(\mathcal{H}_{K}, F_{K})$ 
and the corresponding spectral triple $(\mathcal{H}_{K} , D_{K})$
in future papers.

\section{Fredholm module on $n$-cube}

\subsection{Definition of Fredholm module}
\label{subsec:dfn-n}

In this subsection, we construct a ``good'' Fredholm module on $n$-cubes $\gamma_{n}$. 
For the simplicity, we set 
$\gamma_{n} = [0,e]^{n}$ in $\mathbb{R}^{n}$ with the length of edge $e > 0$; the following construction applies to any $n$-cubes.

Let $V$ be the set of vertices of $\gamma_{n}$:  
\[
V = \{
(a_{1}, \dots , a_{n}) \in \mathbb{R}^{n} \,;\, 
a_{i} = 0 \text{ or } e \quad (i=1,2,\dots , n) \}. 
\]
We give a number of vertices in $V$ inductively. 
For $n=1$, an interval $\gamma_{1} = [0,e]$ has two vertices $0$ and $e$. 
Set $v_{0} = 0$ and $v_{1} = e$. 
For an arbitrary $n$, we assume that we have a number of vertices of $\gamma_{n-1}$. 
Then a number of vertices of $\gamma_{n}$ is as follows: 

\begin{enumerate} 
\item  
$v_{i} = (a_{1}, \dots , a_{n-1}, 0) = (a_{1}, \dots , a_{n-1})$ 
($0 \leq i \leq 2^{n-1}-1$) 
under the inclusion
$\gamma_{n-1} \to \gamma_{n-1} \times \{ 0 \} \subset \gamma_{n}$. 
\item  
$v_{2^{n}-1-i} = (a_{1} , \dots , a_{n-1} , e)$ 
($0 \leq i \leq 2^{n-1}-1$) 
if 
$v_{i} = (a_{1} , \dots , a_{n-1}, 0)$. 
\end{enumerate} 

\begin{example}
\begin{enumerate}
\item 
When $n=2$, the numbering of vertices given by 
$v_{0} = (0,0)$, $v_{1} = (e,0)$, $v_{2} = (e,e)$, $v_{3} = (0,e)$; 
see Figure \ref{fig:n=2}.  
\item 
When $n=3$, the numbering of vertices is given by 
\[
\begin{array}{cccc}
v_{0} = (0,0,0), & v_{1} = (e,0,0), & v_{2} = (e,e,0), & v_{3} = (0,e,0), \\ 
v_{4} = (0,e,e), & v_{5} = (e,e,e), & v_{6} = (e,0,e), & v_{7} = (0,0,e). 
\end{array}
\]
See Figure \ref{fig:n=3}.  
\end{enumerate}
\begin{figure}[h]
\begin{minipage}{0.45\hsize}
\begin{equation*}
\xymatrix{
v_{3} & v_{2} \ar@{-}[l] \\ 
v_{0} \ar@{-}[r] \ar@{-}[u] & v_{1} \ar@{-}[u]
}
\end{equation*}
\caption{$n=2$. }
\label{fig:n=2}
\end{minipage}
\begin{minipage}{0.45\hsize}
\begin{equation*}
\xymatrix{
 & v_{4} \ar@{-}[ld] \ar@{-}[rr] & & v_{5} \ar@{-}[ld] \\ 
v_{7} & & v_{6} \ar@{-}[ll] & \\
 & v_{3} \ar@{.}[uu] & & v_{2} \ar@{.}[ll] \ar@{-}[uu] \\ 
v_{0} \ar@{-}[rr] \ar@{.}[ru] \ar@{-}[uu] & & 
	v_{1} \ar@{-}[ru] \ar@{-}[uu] & 
}
\end{equation*}
\caption{$n=3$.} 
\label{fig:n=3}
\end{minipage}
\end{figure}
\end{example}

Set 
$V_{0} = \{ v_{i} \,;\, i = \text{ even} \}$ and  
$V_{1} = \{ v_{i} \,;\, i = \text{ odd}\}$, 
so we have $V = V_{0} \cup V_{1}$. 
Set also
\begin{align*}
\mathcal{H}^{+} &= \ell^{2}(V_{0}) 
	= \ell^{2}(v_{0}) \oplus \ell^{2}(v_{2}) \oplus \dots \oplus \ell^{2}(v_{2^{n}-2}), \\  
\mathcal{H}^{-} &= \ell^{2}(V_{1}) 
	= \ell^{2}(v_{1}) \oplus \ell^{2}(v_{3}) \oplus \dots \oplus \ell^{2}(v_{2^{n}-1})
\end{align*} 
and  
$\mathcal{H} = \mathcal{H}^{+} \oplus \mathcal{H}^{-}$.  
The vector space 
$\mathcal{H} (\cong \mathbb{C}^{2^{n}})$  is 
a Hilbert space of 
dimension $2^{n}$ with an inner product 
\[
\langle f , g \rangle 
= \sum_{i = 0}^{2^{n}-1} f(v_{i})\overline{g(v_{i})}. 
\]
We assume that $\mathcal{H}$ is $\mathbb{Z}_{2}$-graded 
with the grading $\epsilon = \pm 1$ on $\mathcal{H}^{\pm}$, respectively.  
The $C^{\ast}$-algebra 
$C(V)$  
of continuous functions on $V$   
acts on $\mathcal{H}$ by multiplication:  
\[
\rho(f) = (f(v_{0}) \oplus f(v_{2}) \oplus \dots \oplus f(v_{2^{n}-2})) 
	\oplus (f(v_{1}) \oplus f(v_{3}) \oplus \dots \oplus f(v_{2^{n}-1})).  
\]

A Fredholm operator $F_{n}$ on $\mathcal{H}$ is also defined inductively. 
We set $X_{1} = 1$ 
and 
$X_{n} = \begin{bmatrix} O & X_{n-1} \\ X_{n-1} & O \end{bmatrix} 
\in M_{2^{n-1}}(\mathbb{C})$ ($n \geq 2$). 
We set also $G_{1} = 1$,  
$G_{n} = \begin{bmatrix} 
G_{n-1} & -X_{n-1} \\ X_{n-1} & G_{n-1} 
\end{bmatrix} \in M_{2^{n-1}}(\mathbb{C})$ ($n \geq 2$) 
and  
$U_{n} = \dfrac{1}{\sqrt{n}}G_{n}$ ($n \geq 1$).

\begin{proposition}
\label{prp:unitary_n}
$U_{n}$ is a unitary matrix.  
\end{proposition}

\begin{proof}
Firstly, we have 
\begin{align*}
&\phantom{=} X_{n}G_{n}^{\ast} - G_{n}X_{n} \\ 
&= 
\begin{bmatrix} O & X_{n-1} \\ X_{n-1} & O \end{bmatrix} 
\begin{bmatrix} G_{n-1}^{\ast}& X_{n-1} \\ -X_{n-1} & G_{n-1}^{\ast} \end{bmatrix} 
- 
\begin{bmatrix} G_{n-1} & -X_{n-1} \\ X_{n-1} & G_{n-1} \end{bmatrix} 
\begin{bmatrix} O & X_{n-1} \\ X_{n-1} & O \end{bmatrix} \\ 
&= 
X_{2} \otimes (X_{n-1}G_{n-1}^{\ast} - G_{n-1}X_{n-1} ) \\ 
& = 
\cdots 
= 
X_{n} \otimes (X_{1}G_{1}^{\ast} - G_{1}X_{1})
= O. 
\end{align*}

We prove $U_{n}U_{n}^{\ast} = E_{2^{n}}$ by induction. 
Clearly, $U_{1}= 1$  
is unitary.  
Assume that $U_{n-1}$ is a unitary matrix. 
Then we have 
\begin{align*}
G_{n-1}G_{n-1}^{\ast} + X_{n-1}^{2} 
	&= (n-1)E_{2^{n-2}} + E_{2^{n-2}} = nE_{2^{n-2}}.  
\end{align*}
Thus we obtain 
\begin{align*}
G_{n}G_{n}^{\ast} 
&= \begin{bmatrix} 
G_{n-1} & -X_{n-1} \\ X_{n-1} & G_{n-1} 
\end{bmatrix}
\begin{bmatrix} 
G_{n-1}^{\ast} & X_{n-1} \\ -X_{n-1} & G_{n-1}^{\ast} 
\end{bmatrix}  \\ 
&= \begin{bmatrix} 
G_{n-1}G_{n-1}^{\ast} + X_{n-1}^{2} 
	& G_{n-1}X_{n-1} - X_{n-1}G_{n-1}^{\ast} \\ 
X_{n-1}G_{n-1}^{\ast} - G_{n-1}X_{n-1} 
	& X_{n-1}^{2} + G_{n-1}G_{n-1}^{\ast}
\end{bmatrix}   \\ 
&= \begin{bmatrix} 
G_{n-1}G_{n-1}^{\ast} + X_{n-1}^{2} 
	& (X_{n-1}G_{n-1}^{\ast} - G_{n-1}X_{n-1})^{\ast} \\ 
X_{n-1}G_{n-1}^{\ast} - G_{n-1}X_{n-1} 
	& X_{n-1}^{2} + G_{n-1}G_{n-1}^{\ast}
\end{bmatrix} 
= nE_{2^{n-1}}. 
\end{align*}
Therefore, $U_{n} = \dfrac{1}{\sqrt{n}}G_{n}$ is a unitary matrix. 
\end{proof}

Set $F_{n} = \begin{bmatrix} 
& U_{n}^{\ast} \\ U_{n} & \end{bmatrix} \in M_{2^{n}}(\mathbb{C})$. 
By Proposition \ref{prp:unitary_n}, 
we have 
$F_{n}^{2}=E_{2^{n}}$ and $F_{n}^{\ast} = F_{n}$. 
We consider that $F_{n}$ is a bounded operator on 
a finite dimensional Hilbert space 
\[ 
\mathcal{H} 
= 
(\ell^{2}(v_{0}) \oplus \ell^{2}(v_{2}) \oplus \dots \oplus \ell^{2}(v_{2^{n}-2})) 
\oplus 
(\ell^{2}(v_{1}) \oplus \ell^{2}(v_{3}) \oplus \dots \oplus \ell^{2}(v_{2^{n}-1}))
\cong \mathbb{C}^{2^{n}}
\]
by the left multiplication of a matrix $F_{n}$. 
Because of 
$F_{n}\epsilon + \epsilon F_{n} = O$, 
$(\mathcal{H} , F_{n})$ is an even Fredholm module on 
$C(V)$.

\begin{example}
\begin{enumerate}
\item 
When $n=1$, we have $F_{1} = 
\begin{bmatrix} 0 & 1 \\ 1 & 0 \end{bmatrix}$, 
which is introduced by \cite[Chapter IV. 3. $\varepsilon$]{MR1303779}. 
\item 
When $n=2$, we have 
\[
G_{2}
= \begin{bmatrix} 1 & -1 \\ 1 & 1 \end{bmatrix}, \quad  
U_{2} = \dfrac{1}{\sqrt{2}}
\begin{bmatrix} 1 & -1 \\ 1 & 1 \end{bmatrix}, \quad 
F_{2} = \dfrac{1}{\sqrt{2}}\begin{bmatrix}
 & & 1 & 1 \\ 
 & & -1 & 1 \\ 
1 & -1 & & \\ 
1 & 1 & &  
\end{bmatrix}. 
\]
\item 
When $n=3$, we have 
$G_{3} 
= \begin{bmatrix} 
1 & -1 & 0 & -1 \\ 
1 & 1 & -1 & 0 \\ 
0 & 1 & 1 & -1 \\ 
1 & 0 & 1 & 1 
\end{bmatrix}$, 
$U_{3} = \dfrac{1}{\sqrt{3}}G_{3}$ 
and 
$F_{3} = \begin{bmatrix} 
& U_{3}^{\ast} \\ U_{3} & \end{bmatrix}$. 
\end{enumerate}
\end{example}

\begin{remark}
\label{rmk:G}
The components of $G_{n}$ correspond to 
the following orientation of edges, 
the correspondence is similar to adjacency matrices of oriented graphs. 
When $n=1$,  the orientation of the graph $\gamma_{1} = [0,e]$ 
is from $v_{0} = 0$ to $v_{1}=e$; 
we denote such an orientation by $v_{0} \to v_{1}$. 
Assume that we have the orientation of the edges of $\gamma_{n-1}$. 
\begin{enumerate}
\item 
Assume $0 \leq i,j \leq 2^{n-1}-1$. 
The orientation in $\gamma_{n}$ is from $v_{i}$ to $v_{j}$; $v_{i} \to v_{j}$, when 
the orientation in $\gamma_{n-1}$ is from $v_{i}$ to $v_{j}$. 
Here, we consider that $\gamma_{n-1}$ is a subset in $\gamma_{n}$ under the inclusion 
$\gamma_{n-1} \to \gamma_{n-1} \times \{ 0 \} \subset \gamma_{n}$. 
\item 
$v_{i} \to v_{2^{n}-1-i}$ $(0 \leq i \leq 2^{n-1}-1)$, which means $(a_{1} , \dots , a_{n-1} , 0) \to (a_{1} , \dots , a_{n-1} , e)$. 
\item 
$v_{2^{n}-1-i} \leftarrow v_{2^{n}-1-j}$
if 
$v_{i} \to v_{j}$ $(0 \leq i,j \leq 2^{n-1}-1)$.
\end{enumerate}

\begin{figure}[h]
\begin{minipage}{0.48\hsize}
\begin{equation*}
\xymatrix{
v_{3} & v_{2} \ar@{->}[l] \\ 
v_{0} \ar@{->}[r] \ar@{->}[u]
	& v_{1} \ar@{->}[u]
}
\end{equation*}
\caption{orientation of edges of $\gamma_{2}$}
\end{minipage}
\begin{minipage}{0.48\hsize}
\begin{equation*}
\xymatrix{
 & v_{4} \ar@{->}[ld] \ar@{->}[rr] & & v_{5} \ar@{->}[ld] \\ 
v_{7} & & v_{6} \ar@{->}[ll] & \\
 & v_{3} \ar@{.>}[uu] & & v_{2} \ar@{.>}[ll] \ar@{->}[uu] \\ 
v_{0} \ar@{->}[rr] \ar@{.>}[ru] \ar@{->}[uu] & & 
	v_{1} \ar@{->}[ru] \ar@{->}[uu] & 
}
\end{equation*}
\caption{orientation of edges of $\gamma_{3}$}
\end{minipage}
\end{figure}

Then the $(i,j)$-component $g_{ij}$ $(1 \leq i,j \leq 2^{n-1})$ of $G_{n}$ 
is as follows. 
\begin{enumerate}
\item 
$g_{ij} = 1$ when $v_{2j-2} \to v_{2i-1}$. 
\item 
$g_{ij} = -1$ when   
$v_{2j-2} \leftarrow  v_{2i-1}$. 
\item 
$g_{ij} = 0$ when 
$v_{2j-2}$ and $v_{2i-1}$ do not connect by an edge. 
\end{enumerate}
\end{remark}

\subsection{Calculation of quantized differential form}
\label{subsec:cal-n}

In this subsection we calculate an operator 
$[F_{n} , x^{\alpha}]$ for the coordinate function $x^{\alpha}$ 
on $\mathbb{R}^{n}$ ($\alpha = 1, 2, \dots , n$). 
We also show  they satisfy a relation of the Clifford algebra 
on the Euclidean vector space of dimension $n$.

Set 
$d_{n}f = [F_{n}, f] 
= \begin{bmatrix} & d_{n}^{-}f \\ d_{n}^{+}f & \end{bmatrix}$, 
and we have 
\begin{align*}
d_{n}^{+}f &= Uf^{+} - f^{-}U \\ 
d_{n}^{-}f &= U^{\ast}f^{-} - f^{+}U^{\ast} 
= - (U\bar{f}^{+} - \bar{f}^{-}U)^{\ast}
= - {}^{t}d_{n}^{+}f, 
\end{align*}
where 
$f^{+} = f|_{V_{0}}$ and $f^{-} = f|_{V_{1}}$. 
Denote by $A \circ B = [a_{ij}b_{ij}]$ 
the Hadamard product of two matrices 
$A = [a_{ij}]$ and $B = [b_{ij}]$ of the same size.

\begin{proposition}
\label{prp:df}
For any $f \in C(V)$,   
we set $f_{a,b} = f(v_{a}) - f(v_{b})$ and 
\[
\Delta_{n}f = [f_{2j,2i+1}]_{i,j=0,1,\dots , 2^{n-1}-1} 
\in \mathcal{B}(\ell^{2}(V_{0}) , \ell^{2}(V_{1})) \cong M_{2^{n-1}}(\mathbb{C}). 
\]
We have 
\[
d_{n}f = \frac{1}{\sqrt{n}}
\begin{bmatrix} 
 & -{}^{t}(\Delta_{n}f \circ G_{n}) \\ 
\Delta_{n}f \circ G_{n} & 
\end{bmatrix}. 
\]
\end{proposition}

\begin{proof}  
As in Remark \ref{rmk:G}, we denote 
$G_{n} = [g_{ij}]$. We have 
\begin{align*}
&\phantom{=} \sqrt{n} d_{n}^{+}f \\ 
&= 
G_{n}
\begin{bmatrix} 
f(v_{0}) & & & \\ 
& f(v_{2}) & & \\ 
& & \ddots & \\ 
& & & f(v_{2^{n}-2})
\end{bmatrix} 
- 
\begin{bmatrix} 
f(v_{1}) & & & \\ 
& f(v_{3}) & & \\ 
& & \ddots & \\ 
& & & f(v_{2^{n}-1})
\end{bmatrix} 
G_{n} \\ 
&= 
\begin{bmatrix} g_{ij}f(v_{2j}) \end{bmatrix} 
- \begin{bmatrix} f(v_{2i-1})g_{ij} \end{bmatrix} \\ 
&= \begin{bmatrix} f_{2j, 2i-1}g_{ij} \end{bmatrix} \\ 
&= \Delta_{n}f \circ G_{n}.  
\end{align*}
\end{proof} 

Thus an $(i,j)$-component of $d_{n}^{+}f$ is $0$ if 
$v_{2i-1}$ and $v_{2j}$ do not connect by an edge. 

\begin{proposition}
\label{prp:e_a}
For the coordinate function 
$x^{\alpha}$ on $\mathbb{R}^{n}$ ($\alpha = 1,2, \dots , n$), 
we set $\displaystyle e_{(n)}^{\alpha} = \frac{\sqrt{n}}{e}d_{n}x^{\alpha}$.  
We have 
\begin{equation}
\label{eq:e_a}
e_{(n)}^{\alpha} = \frac{\sqrt{n}}{e}d_{n}x^{\alpha} 
= 
\begin{bmatrix} 
 & E_{2^{n-\alpha-1}} \otimes \begin{bmatrix} 1 & \\ & -1 \end{bmatrix} \\ 
-E_{2^{n-\alpha-1}} \otimes \begin{bmatrix} 1 & \\ & -1 \end{bmatrix} &  
\end{bmatrix} 
\otimes X_{\alpha}. 
\end{equation}
Here,  
$E_{1/2} \otimes \begin{bmatrix} 1 & \\ & -1 \end{bmatrix} = 1$. 
\end{proposition}

\begin{proof} 
Firstly, by using 
$\Delta_{n}x^{n} \circ G_{n} = -e X_{n}$ 
and Proposition \ref{prp:df}, we have 
$e_{(n)}^{n} = 
\begin{bmatrix} & X_{n} \\ -X_{n} & \end{bmatrix}$. 

Next we calculate 
$e_{(n)}^{n-1+} =  \dfrac{\sqrt{n}}{e}d_{n}^{+}x^{n-1}$. 
By the definition of the numbering of vertices and 
the orientation of edges of $\gamma_{n}$, 
for $0 \leq i,j \leq 2^{n}-1$,  
``$v_{i} \to v_{j}$ is positive (resp. negative) with $x^{n-1}$ direction''  
if and only if 
``$v_{i + 2^{n-1}} \to v_{j + 2^{n-1}}$ is 
negative (resp. positive) with $x^{n-1}$ direction''. 
So we have 
$e_{(n)}^{n-1+} = 
\begin{bmatrix} e_{(n-1)}^{n-1+} & \\ & -e_{(n-1)}^{n-1+} \end{bmatrix}
=
\begin{bmatrix} 1 & \\ & -1 \end{bmatrix} \otimes (-X_{n-1})$. 
This implies 
\[
e_{(n)}^{n-1} = 
\begin{bmatrix} 
 & E_{1} \otimes \begin{bmatrix} 1 & \\ & -1 \end{bmatrix} \\ 
-E_{1} \otimes \begin{bmatrix} 1 & \\ & -1 \end{bmatrix} &  
\end{bmatrix} \otimes X_{n-1}. 
\]

We calculate $e_{(n)}^{\alpha}$ ($\alpha = 1,2, \dots , n-2$) by induction on $n \geq 3$. 
Note that the calculation of $e_{(n)}^{\alpha}$ for $n=1, 2$ is already done. Namely, the beginning of induction is the following: 
\[
e_{(1)}^{1} = \begin{bmatrix} & 1 \\ -1 & \end{bmatrix}, \quad   
e_{(2)}^{1} = 
\begin{bmatrix} 
 & & 1 & \\ 
 & & & -1 \\ 
-1 & & & \\ 
 & 1 & & 
\end{bmatrix}, \quad 
e_{(2)}^{2} = 
\begin{bmatrix} 
 & & & 1 \\ 
 & & 1 & \\ 
 & -1 & & \\ 
-1 & & & 
\end{bmatrix}. 
\]

Assume that  equation (\ref{eq:e_a}) holds for $n-1$.   
By the definition of the numbering of vertices and 
the orientation of edges of $\gamma_{n}$, 
for $1 \leq \alpha \leq n-2$, 
``$i \to j$ is positive (resp. negative) with $x^{\alpha}$ direction''  
if and only if 
``$v_{i + 2^{n-1}} \rightarrow v_{j + 2^{n-1}}$ is 
positive (resp. negative) with $x^{\alpha}$ direction''. 
So we have 
\begin{align*} 
e_{(n)}^{\alpha+} 
&= 
\begin{bmatrix} e_{(n-1)}^{\alpha+} & \\ & e_{(n-1)}^{\alpha+} \end{bmatrix}
= E_{2} \otimes e_{(n-1)}^{\alpha+}  
= -E_{2} \otimes 
\left( E_{2^{n-1-\alpha-1}} \otimes \begin{bmatrix} 1 & \\ & -1 \end{bmatrix} \otimes X_{\alpha} \right)  \\ 
&= - E_{2^{n-\alpha-1}} \otimes \begin{bmatrix} 1 & \\ & -1 \end{bmatrix} \otimes X_{\alpha}.  
\end{align*} 
Therefore we have 
\[
e_{(n)}^{\alpha} = 
\begin{bmatrix} 
 & E_{2^{n-\alpha-1}} \otimes \begin{bmatrix} 1 & \\ & -1 \end{bmatrix} \\ 
-E_{2^{n-\alpha-1}} \otimes \begin{bmatrix} 1 & \\ & -1 \end{bmatrix} &  
\end{bmatrix} \otimes X_{\alpha} 
\quad (\alpha = 1,2,\dots , n-2). 
\]

We have equation (\ref{eq:e_a}) by the above calculations for 
any $n$ and $\alpha = 1,2,\dots n$.  
\end{proof}

By the explicit formula of $e_{(n)}^{\alpha}$ 
in Proposition \ref{prp:e_a}, 
we have a Clifford relation of $d_{n}x^{\alpha}$. 

\begin{proposition}
\label{prp:Clifford}
We have 
\[
e_{(n)}^{\alpha}e_{(n)}^{\beta} 
= 
\begin{cases}
-e_{(n)}^{\beta}e_{(n)}^{\alpha} & (\alpha \neq \beta) \\ 
-E_{2^{n}} & (\alpha = \beta) 
\end{cases}. 
\]
By $d_{n}x^{\alpha} = \dfrac{e}{\sqrt{n}}e_{(n)}^{\alpha}$,  
we have 
\[
d_{n}x^{\alpha}d_{n}x^{\beta}
= 
\begin{cases}
-d_{n}x^{\beta}d_{n}x^{\alpha} & (\alpha \neq \beta) \\ 
-\dfrac{e^{2}}{n}E_{2^{n}} & (\alpha = \beta) 
\end{cases}. 
\]
\end{proposition} 

\begin{proof} 
Firstly, we have  
\[ 
e_{(n)}^{\alpha}e_{(n)}^{\alpha} 
= 
\begin{bmatrix} 
-E_{2^{n-\alpha-1}}^{2} \otimes \begin{bmatrix} 1 & \\ & -1 \end{bmatrix}^{2} &  \\ 
& -E_{2^{n-\alpha-1}}^{2} \otimes \begin{bmatrix} 1 & \\ & -1 \end{bmatrix}^{2}   
\end{bmatrix} 
\otimes X_{\alpha}^{2}
= -E_{2^{n}}. 
\]  

Set $k = \alpha - \beta > 0$, then we have 
$X_{\alpha} = X_{k+1} \otimes X_{\beta}$. 
We can rewrite $e_{(n)}^{\alpha}$ and $e_{(n)}^{\beta}$ as follows: 
\begin{align*}
e_{(n)}^{\alpha}
&= 
\begin{bmatrix} 
 & E_{2^{n-\alpha-1}} \otimes \begin{bmatrix} X_{k+1} & \\ & -X_{k+1} \end{bmatrix} \\ 
-E_{2^{n-\alpha-1}} \otimes \begin{bmatrix} X_{k+1} & \\ & -X_{k+1} \end{bmatrix} &  
\end{bmatrix} \otimes X_{\beta},  \\ 
e_{(n)}^{\beta} 
&= 
\begin{bmatrix} 
 & E_{2^{n-\alpha-1}} \otimes E_{2^{k}} \otimes \begin{bmatrix} 1 & \\ & -1 \end{bmatrix} \\ 
-E_{2^{n-\alpha-1}} \otimes E_{2^{k}} \otimes \begin{bmatrix} 1 & \\ & -1 \end{bmatrix} &  
\end{bmatrix} \otimes X_{\beta}. 
\end{align*} 
Now, we set 
$\epsilon_{1} = \begin{bmatrix} 1 & \\ & -1 \end{bmatrix}$, and we have 
\begin{align*}
\begin{bmatrix} X_{k+1} & \\ & -X_{k+1} \end{bmatrix} 
\left(E_{2^{k}} \otimes \begin{bmatrix} 1 & \\ & -1 \end{bmatrix}\right) 
&= 
\begin{bmatrix} 
X_{k} \otimes (X_{2}\epsilon_{1}) & \\ 
 & -X_{k} \otimes (X_{2} \epsilon_{1}) 
\end{bmatrix} \text{ and}  \\ 
\left(E_{2^{k}} \otimes \begin{bmatrix} 1 & \\ & -1 \end{bmatrix}\right) 
\begin{bmatrix} X_{k+1} & \\ & -X_{k+1} \end{bmatrix} 
&= 
\begin{bmatrix} 
X_{k} \otimes (\epsilon_{1} X_{2}) & \\ 
 & -X_{k} \otimes (\epsilon_{1} X_{2}) 
\end{bmatrix}. 
\end{align*}
Thus the relation
$e_{(n)}^{\alpha}e_{(n)}^{\alpha} = -e_{(n)}^{\alpha}e_{(n)}^{\alpha}$ ($\alpha \neq \beta$) 
holds
since we have 
$X_{2}\epsilon_{1} + \epsilon_{1} X_{2} = O$. 
\end{proof}

\begin{remark}
When we take the limit as the length of edges tends to $0$, that is $e \to 0$, 
we have 
\[
d_{n}x^{\alpha}d_{n}x^{\alpha} = -\dfrac{e^{2}}{n}E_{2^{n}} \to O. 
\]
Thus we regard $d_{n}x^{\alpha}$ as 
a quantization of the ordinal exterior differential $dx^{\alpha}$ on $\mathbb{R}^{n}$. 
\end{remark}

\begin{remark}
\label{rmk:unitary}
For any unitary matrix $U \in U(2^{n-1})$, 
since an odd matrix 
$F = \begin{bmatrix} & U^{\ast} \\ U & \end{bmatrix}$  
defines an operator on $\mathcal{H}$, $F$ defines a Fredholm module on $C(V)$. 
Moreover, 
since any $F$ is homotopic to $F_{n}$,  
it defines a same $K$-homology class in $K^{0}(C(V))$. 
However, the general $F$ sometimes does not have good properties. 
For example, we have  $[F , x^{\alpha}] = O$ for $\alpha = 2, 3, \dots , n$ when we assume $U = E_{2^{n-1}}$ . 
Thus, in this case, we cannot regard $[F, x^{\alpha}]$ as 
a quantization of the ordinal exterior differential $dx^{\alpha}$ on $\mathbb{R}^{n}$. 
\end{remark}

By Proposition \ref{prp:Clifford}, 
we have the volume element 
$\omega_{n} = e_{(n)}^{1}e_{(n)}^{2} \cdots e_{(n)}^{n}$ 
in the Clifford algebra.  
We can easily calculate its absolute value $|\omega_{n}|$. 
We do not use $|\omega_{n}|$ directly, but we use  
$|d_{n}x^{1} d_{n}x^{2} \cdots d_{n}x^{n}|$, 
which is a constant multiple of $|\omega_{n}|$; see also Section \ref{subsec:Dixqv}.

\begin{proposition}
\label{prp:nDixmier}
We have
$
|[F_{n} , x^{1}] \cdots [F_{n} , x^{n}]| = 
\dfrac{e^{n}}{n^{n/2}}E_{2^{n}}
$. 
By the definition of $e_{(n)}^{\alpha}$, we also have $|\omega_{n}| = E_{2^{n}}$. 
\end{proposition}

\begin{proof}  
Because of  
$[F_{n}, x^{\alpha}]^{\ast}[F_{n}, x^{\alpha}] 
= \dfrac{e^{2}}{n}e_{(n)}^{\alpha \ast}e_{(n)}^{\alpha} 
= \dfrac{-e^{2}}{n}(e_{(n)}^{\alpha})^{2} 
= \dfrac{e^{2}}{n}E_{2^{n}}$,  
we have 
\begin{align*}
|[F_{n} , x^{1}] \cdots [F_{n} , x^{n}]|^{2} 
&= 
([F_{n} , x^{1}] \cdots [F_{n} , x^{n}])^{\ast}[F_{n} , x^{1}] \cdots [F_{n} , x^{n}] \\ 
&= 
[F_{n} , x^{n}]^{\ast} \cdots [F_{n} , x^{1}]^{\ast}  
[F_{n} , x^{1}] \cdots [F_{n} , x^{n}] \\ 
&= 
\left( \dfrac{e^{2}}{n} \right)^{n} E_{2^{n}}. 
\end{align*}
This implies 
\begin{align*}
|[F_{n} , x^{1}] \cdots [F_{n} , x^{n}]| = 
\frac{e^{n}}{n^{n/2}}E_{2^{n}}. 
\end{align*}

\end{proof}

\section{Fredholm module on self-similar sets built on $n$-cubes}
\label{sec:FonK}

\subsection{Fredholm module and spectral triple}
\label{subsec:FMST}

In this subsection, we construct a Fredholm module and a spectral triple on 
self-similar sets built on any $n$-cubes $\gamma_{n}$. 
For the simplicity, we assume that the length of edges of $\gamma_{n}$ equals $1$.  
Let 
$f_{s} : \gamma_{n} \to \gamma_{n}$ 
($s=1, \dots N$) 
be similitudes.  
We define the similarity ratio of $f_{s}$ to be
\[
r_{s} = \dfrac{\|f_{s}(x) - f_{s}(y)\|_{\mathbb{R}^{n}}}{\| x - y \|_{\mathbb{R}^{n}}} \;\,
(<1)
\quad (x \neq y).
\]  
An iterated function system (IFS)
$(\gamma_{n}, S = \{ 1, \dots , N \}, \{ f_{s} \}_{s \in S})$ 
defines the unique non-empty compact set 
$K = K(\gamma_{n}, S = \{ 1, \dots , N \}, \{ f_{s} \}_{s \in S})$ 
called the self-similar set 
such that 
$K = \bigcup_{s=1}^{N} f_{s}(K)$.  
We use $\dim_{S}(K)$ to denote the similarity dimension of $K$, that is, the number $s$ that satisfies 
\[
\sum_{s=1}^{N}r_{s}^{s} = 1. 
\]
If an IFS
$(\gamma_{n}, S, \{ f_{s} \}_{s \in S})$ 
satisfies 
the open set condition, $\dim_{S}(K)$ turns out to be equal to the Hausdorff dimension $\dim_{H}(K)$ of $K$.

Set $f_{\bs{s}} = f_{s_{1}} \circ \dots \circ f_{s_{j}}$ for 
$\bs{s} = (s_{1}, \dots ,  s_{j}) \in S^{\infty} = \bigcup_{j=0}^{\infty} S^{\times j}$ and 
$f_{\emptyset} = \mathrm{id}$. 
For the simplicity, 
we will use $i$ to express
the vertex $f_{\bs{s}} (v_{i})$ 
of 
an $n$-cube $f_{\bs{s}} (\gamma_{n})$ and write $V_{\bs{s}}$ as the vertices of an $n$-cube $f_{\bs{s}} (\gamma_{n})$. 
We also denote the length of the edge of $f_{\bs{s}} (\gamma_{n})$ by $e_{\bs{s}}$. 
As introduced in subsection \ref{subsec:dfn-n}, 
we define the Hilbert space $\mathcal{H}_{\bs{s}} = \ell^{2}(V_{\bs{s}})$ 
on an $n$-cube of the length $e_{\bs{s}}$ that consists of the positive part
$\mathcal{H}_{\bs{s}}^{+}$ 
and the negative part 
$\mathcal{H}_{\bs{s}}^{-}$.  
By taking the direct sum on all $n$-cubes, we define the following data: 
\[ 
\mathcal{H}_{K} = \bigoplus_{\bs{s} \in S^{\infty}}\mathcal{H}_{\bs{s}}, 
\quad 
F_{K} = \bigoplus_{\bs{s} \in S^{\infty}} F_{n}, \quad   
D_{K} = \bigoplus_{\bs{s} \in S^{\infty}} \dfrac{1}{e_{\bs{s}}}F_{n}. 
\]

Let $V_{K}$ be the closure of the set of vertices of all $n$-cubes 
$f_{\bs{s}} (\gamma_{n}) \subset \mathbb{R}^{n} $. 
That is, $V_{K}$ is the closure of 
$\displaystyle  
\bigcup_{\bs{s} \in S^{\infty}}
V_{\bs{s}}$. 
Then, if $\displaystyle V \subset \bigcup_{s=1}^{N}f_{s}(V)$ holds, 
we have $V_{K} = K$. 
If not, $V_{K}$ equals the union of 
$\displaystyle 
\bigcup_{\bs{s} \in S^{\infty}}
V_{\bs{s}}$ and $K$. 
We also let 
$\mathcal{A}_{K}$   
be the Banach algebra of Lipschitz functions $\mathrm{Lip}(V_{K})$ 
on $V_{K}$ with the norm 
$\| a \|_{\mathcal{A}_{K}} = \| a \|_{\infty} + \mathrm{Lip}(a)$, 
where the second term is the Lipschitz constant of a Lipschitz function $a$.  
The Banach algebra $\mathcal{A}_{K}$ acts on  
$\mathcal{H}_{K}$ by 
\[
\rho_{K} : \mathcal{A}_{K} \to \mathcal{B}(\mathcal{H}_{K}) \,;\, 
\rho_{K}(a)(\oplus \xi_{\bs{s}} )  = \oplus (a|_{V_{\bs{s}}}) \cdot \xi_{\bs{s}} . 
\] 

\begin{lemma}
\label{lem:KDirac}
Define
\[ \mathcal{H}_{K}^{1}= \left\{
\bigoplus_{\bs{s} \in S^{\infty}} \xi_{\bs{s}} \in \mathcal{H}_{K}  \,;\, 
\| \oplus \xi_{\bs{s}}  \|_{\mathcal{H}_{K}^{1}}^{2} 
= 
\sum_{\bs{s} \in S^{\infty}} \frac{1}{e_{\bs{s}}^{2}}\sum_{i=0}^{2^{n}-1}|\xi_{\bs{s}}  (i)|^{2} 
< \infty 
\right\}.
\]
Then, an operator 
$D_{K}$ is a self-adjoint operator of $\text{dom}(D_{K}) = \mathcal{H}_{K}^{1}$.  
\end{lemma}

\begin{proof} 
By inclusions 
$\{ \oplus \xi_{\bs{s}}  \in \mathcal{H}_{K} \,;\,  \xi_{\bs{s}}  = 0  \text{ except finite } \bs{s} \} 
\subset \mathcal{H}_{K}^{1} \subset \mathcal{H}_{K}$, 
$\mathcal{H}_{K}^{1}$ is a dense subset in $\mathcal{H}_{K}$. 

On each $n$-cubes $f_{\bs{s}} (\gamma_{n})$,  
we have 
\begin{equation*}
\| F_{n} \xi_{\bs{s}}  \|_{\ell^{2}}^{2} 
= 
\| U_{n}\xi_{\bs{s}} ^{+} \|_{\ell^{2}}^{2} + \| U_{n}^{\ast}\xi_{\bs{s}} ^{-} \|_{\ell^{2}}^{2} 
= 
\| \xi_{\bs{s}} ^{+} \|_{\ell^{2}}^{2} + \| \xi_{\bs{s}} ^{-} \|_{\ell^{2}}^{2} 
= 
\sum_{i=0}^{2^{n}-1} |\xi_{\bs{s}} (i)|^{2}
\end{equation*} 
for any function $\xi_{\bs{s}} $ on 
$V_{\bs{s}}$, 
where $\xi_{\bs{s}} ^{\pm}$ denote the  
$\mathcal{H}_{\bs{s}}^{\pm}$ parts of $\xi_{\bs{s}} $, respectively.  
Then, we have 
\[ 
\| D_{K} (\oplus \xi_{\bs{s}} ) \|_{\mathcal{H}_{K}}^{2} =  
\sum_{\bs{s} \in S^{\infty}} \frac{1}{e_{\bs{s}}^{2}}\sum_{i=0}^{2^{n}-1}|\xi_{\bs{s}} (i)|^{2} 
= 
\| \oplus \xi_{\bs{s}}  \|_{\mathcal{H}_{K}^{1}}^{2}
\]
for $\oplus \xi_{\bs{s}}  \in \mathcal{H}_{K}$. 
Thus we have 
$D_{K} (\mathcal{H}_{K}^{1}) \subset \mathcal{H}_{K}$, and $D_{K}$ is a symmetric operator with domain $\mathcal{H}_{K}^{1}$. 

On the other hand,  
we set 
$\oplus \eta_{\bs{s}}  = \oplus e_{\bs{s}}F_{n}\xi_{\bs{s}} $ 
for any $\oplus \xi_{\bs{s}}  \in \mathcal{H}_{K}$. Then, 
$\oplus \eta_{\bs{s}}  \in \mathcal{H}_{K}^{1}$ since  
\begin{align*}
\| \oplus \eta_{\bs{s}}  \|_{\mathcal{H}_{K}^{1}}^{2} 
&= 
\sum_{\bs{s} \in S^{\infty}}  \|F_{n}\xi_{\bs{s}}\|_{\ell^{2}}^{2} 
= 
\sum_{\bs{s} \in S^{\infty}} \|\xi_{\bs{s}}\|_{\ell^{2}}^{2} 
= 
\| \oplus \xi_{\bs{s}}  \|_{\mathcal{H}_{K}}^{2} < \infty  . 
\end{align*} 
This implies 
$D_{K} (\mathcal{H}_{K}^{1}) \supset \mathcal{H}_{K}$. 
Thus we have 
$D_{K} (\mathcal{H}_{K}^{1}) = \mathcal{H}_{K}$. 
Therefore $D_{K}$ is a self-adjoint operator of 
\[
\text{dom}(D_{K}) = 
 \left\{
\bigoplus_{\bs{s} \in S^{\infty}} \xi_{\bs{s}} \in \mathcal{H}_{K}  \,;\, 
\sum_{\bs{s} \in S^{\infty}} \frac{1}{e_{\bs{s}}^{2}}\sum_{i=0}^{2^{n}-1}|\xi_{\bs{s}}(i)|^{2}
< \infty 
\right\} . 
\]

\end{proof} 

Note that we have 
$\rho_{K}(\mathcal{A}_{K})(\mathcal{H}_{K}^{1}) \subset \mathcal{H}_{K}^{1}$ 
and 
$F_{K} = D_{K}|D_{K}|^{-1}$. 
We now prove some regularity of $F_{K}$ and $D_{K}$.

\begin{lemma}
\label{lem:Ksummable}
We have the followings: 
\begin{enumerate} 
\item 
$[F_{K}  , a] \in \mathcal{K}(\mathcal{H}_{K})$ for any $a \in C(V_{K})$. 
\item 
$[D_{K} , a] \in \mathcal{B}(\mathcal{H}_{K})$ for any $a \in \mathcal{A}_{K}$. 
\item 
$|D_{K}|^{-1} \in \mathcal{K}(\mathcal{H}_{K})$. 
\item 
$(D_{K}^{2}+1)^{-1/2} \in \mathcal{K}(\mathcal{H}_{K})$. 
\item 
$|D_{K}|^{-p} \in \mathcal{L}^{1}(\mathcal{H}_{K}) \iff 
p > \dim_{S}(K)$, where 
$\mathcal{L}^{1}(\mathcal{H}_{K})$ is the set of trace class operators 
on $\mathcal{H}_{K}$.  
\item 
$(D_{K}^{2}+1)^{-p/2} \in \mathcal{L}^{1}(\mathcal{H}_{K}) \iff 
	p > \dim_{S}(K)$.  
\end{enumerate}
\end{lemma}

\begin{proof} 
\begin{enumerate}
\item 
First, we take $a \in \mathcal{A}_{K}$. 
For any $\bs{s} \in S^{\times j}$, 
we have 
\begin{align*}
[F_{K} , a ]|_{\mathcal{H}_{\bs{s}}} 
= 
\frac{1}{\sqrt{n}} 
\begin{bmatrix} 
 & -{}^{t}(\Delta_{n}a \circ G_{n}) \\ 
\Delta_{n}a \circ G_{n} & 
\end{bmatrix} . 
\end{align*}
Therefore, the operator norm $\| [F_{K} , a ]|_{\mathcal{H}_{\bs{s}}}\|$ 
is less than 
\begin{align*}
\mathrm{Lip}(a) \cdot e_{\bs{s}} 
= \mathrm{Lip}(a) \cdot \prod_{k=1}^{j}r_{s_{k}}. 
\end{align*}
Thus $[F_{K} , a]$ is compact for $a \in \mathcal{A}_{K}$ 
since we have 
$\displaystyle \prod_{k=1}^{j}r_{s_{k}} 
\leq \max_{s \in S} r_{s}^{j} \to 0$ as $j \to \infty$. 
The case for any continuous function is proved by the denseness of $\mathcal{A}_{K}$ in $C(V_{K})$. 
\item 
For any $\bs{s} \in S^{\times j}$, 
we have 
\begin{align*}
[D_{K} , a ]|_{\mathcal{H}_{\bs{s}}} 
= 
\frac{1}{\sqrt{n}} 
\left( \prod_{k=1}^{j}r_{s_{k}} \right)^{-1}
\begin{bmatrix} 
 & -{}^{t}(\Delta_{n}a \circ G_{n}) \\ 
\Delta_{n}a \circ G_{n} & 
\end{bmatrix}.  
\end{align*}
So the  operator norm $\| [D_{K} , a ]|_{\mathcal{H}_{\bs{s}}}\|$ 
is less than $
\mathrm{Lip}(a)$,  
which is independent of $j$. 
Therefore $[D_{K} , a]$
is bounded on $\mathcal{H}_{K}$. 
\item 
Because of  
$\displaystyle 
|D_{K}| = \bigoplus_{\bs{s} \in S^{\infty}} 
\frac{1}{e_{\bs{s}}}E_{2^{n}}$, 
we have  
$\displaystyle 
|D_{K}|^{-1} = \bigoplus_{j=0}^{\infty} \bigoplus_{\bs{s} \in S^{\times j}}
\left( \prod_{k=1}^{j}r_{s_{k}}\right)E_{2^{n}}$. 
Thus 
$|D_{K}|^{-1}$ is compact since 
we have 
$\displaystyle  \prod_{k=1}^{j}r_{s_{k}} \to 0$ as $j \to \infty$. 
\item 
Because of  
$\displaystyle 
D_{K}^{2}+1 = \bigoplus_{\bs{s} \in S^{\infty}} 
\left(\frac{1}{e_{\bs{s}}^{2}} + 1\right)E_{2^{n}}$, 
we have 
\begin{align*}
(D_{K}^{2}+1)^{-1/2} 
&= 
\bigoplus_{j=0}^{\infty} 
\bigoplus_{\bs{s} \in S^{\times j}}
\left(\prod_{k=1}^{j}r_{s_{k}}^{-2} +1\right)^{-1/2}E_{2^{n}}.   
\end{align*}
Thus 
$(D_{K}^{2}+1)^{-1/2}$ 
is a compact operator. 
\item 
Because of
$\displaystyle 
|D_{K}|^{-p} = \bigoplus_{j=0}^{\infty} 
\bigoplus_{\bs{s} \in S^{\times j}}
\left( \prod_{k=1}^{j}r_{s_{k}}^{p}\right)E_{2^{n}}$,  
we have 
\begin{align*}
\mathrm{Tr}(|D_{K}|^{-p}) 
&= 
\sum_{j=0}^{\infty}\sum_{\bs{s} \in S^{\times j}} 
	2^{n} \prod_{k=1}^{j}r_{s_{k}}^{p}
= 
2^{n} \sum_{j=0}^{\infty} \left( \sum_{s=1}^{N}r_{s}^{p}\right)^{j}. 
\end{align*}
Thus we have 
\begin{equation*}
|D_{K}|^{-p} \in \mathcal{L}^{1}(\mathcal{H}_{K}) 
\iff 
\sum_{s=1}^{N}r_{s}^{p}   < 1. 
\end{equation*}
This implies 
$
|D_{K}|^{-p} \in \mathcal{L}^{1}(\mathcal{H}_{K})
\iff 
p > \dim_{S}(K) $. 
\item 
Because of 
\begin{align*}
(D_{K}^{2}+1)^{-p/2} 
&= 
\bigoplus_{j=0}^{\infty} 
\bigoplus_{\bs{s} \in S^{\times j}}
\left(\prod_{k=1}^{j}r_{s_{k}}^{-2} +1\right)^{-p/2}E_{2^{n}},  
\end{align*}
we have 
\begin{align*}
\mathrm{Tr}((D_{K}^{2}+1)^{-p/2}) 
&= 
\sum_{j=0}^{\infty}\sum_{\bs{s} \in S^{\times j}} 
	2^{n}\left(\prod_{k=1}^{j}r_{s_{k}}^{-2} +1\right)^{-p/2}. 
\end{align*}
Thus we have 
\begin{align*}
\sum_{j=0}^{\infty}\sum_{\bs{s} \in S^{\times j}} 
	2^{n-p/2}\prod_{k=1}^{j}r_{s_{k}}^{p}
\leq 
\mathrm{Tr}((D_{K}^{2}+1)^{-p/2}) 
\leq 
\sum_{j=0}^{\infty}\sum_{\bs{s} \in S^{\times j}} 
	2^{n} \prod_{k=1}^{j}r_{s_{k}}^{p}, 
\end{align*}
that is we have 
\begin{align*}
2^{n-p/2}
\sum_{j=0}^{\infty}\left( \sum_{s=1}^{N} r_{s_{k}}^{p}\right)^{j}
\leq 
\mathrm{Tr}((D_{K}^{2}+1)^{-p/2}) 
\leq 
2^{n}
\sum_{j=0}^{\infty}\left( \sum_{s=1}^{N} r_{s_{k}}^{p}\right)^{j}. 
\end{align*}
This implies 
\begin{align*}
(D_{K}^{2}+1)^{-p/2} \in \mathcal{L}^{1}(\mathcal{H}_{K}) 
\iff 
\sum_{s=1}^{N} r_{s_{k}}^{p} < 1 
\iff 
p > \dim_{S}(K). 
\end{align*}
\end{enumerate}

\end{proof} 

\begin{theorem}
\label{thm:Fmodule} 
The pair 
$(\mathcal{H}_{K}, F_{K})$ is an  
even Fredholm module over $C(V_{K})$ with the $\mathbb{Z}_{2}$-grading 
$\epsilon_{K} = \bigoplus_{\bs{s} \in S^{\infty}} \epsilon$. 
The pair 
$(\mathcal{H}_{K}, F_{K})$ is a 
$([\dim_{S}(K)]+1)$-summable even Fredholm module 
over $\mathcal{A}_{K}$. 
In particular, if we have 
$\dim_{S}(K) < n$, 
an operator 
\[
[F_{K} , a^{1}][F_{K} , a^{2}]  \cdots [F_{K} , a^{n}] 
\]
is of trace class for any 
$a^{1}, a^{2}, \dots , a^{n} \in \mathcal{A}_{K}$. 
\end{theorem}

\begin{proof} 
By the definition of $F_{K}$, we have 
$F_{K}^{2} =1$, 
$F_{K}^{\ast} = F_{K}$ and 
$F_{K}\epsilon_{K} + \epsilon_{K}F_{K} = 0$. 
$[F_{K}, a]$ is also a compact operator 
by Lemma \ref{lem:Ksummable}. 
Therefore, 
$(\mathcal{H}_{K}, F_{K})$
is an even Fredholm module over $C(V_{K})$. 

Next we prove summability of the Fredholm module 
$(\mathcal{H}_{K}, F_{K})$ over $\mathcal{A}_{K}$. 
Since 
$[D_{K} , a]$ is a bounded operator for $a \in \mathcal{A}_{K}$ and 
$|D_{K}|^{-([\dim_{S}(K)]+1)}$ is of trace class, 
we have 
\begin{align*}
& [F_{K} , a^{1}][F_{K} , a^{2}]  \cdots [F_{K} , a^{[\dim_{S}(K)]+1}]  \\ 
=& 
[D_{K} , a^{1}]|D_{K}|^{-1} 
[D_{K}, a^{2}]|D_{K}|^{-1} 
  \cdots [D_{K} , a^{[\dim_{S}(K)]+1}]|D_{K}|^{-1} \\ 
=& 
[D_{K} , a^{1}][D_{K}, a^{2}] \cdots  [D_{K} , a^{[\dim_{S}(K)]+1}]
|D_{K}|^{-([\dim_{S}(K)]+1)} 
\in \mathcal{L}^{1}(\mathcal{H}_{K})
\end{align*} 
for 
$a^{1}, a^{2}, \dots , a^{[\dim_{S}(K)]+1} \in \mathcal{A}_{K}$.  
Here, we have 
$[|D_{K}|^{-1} , T] = 0$ if $T \in \mathcal{B}(\mathcal{H}_{K})$ 
is a direct sum of operators on each $n$-cubes $f_{\bs{s}}(\gamma_{n})$. 
Therefore, 
$(\mathcal{H}_{K}, F_{K})$ is a 
$([\dim_{S}(K)]+1)$-summable even Fredholm module.  

\end{proof} 

\begin{theorem}
\label{thm:triple} 
The tirple  
$(\mathcal{A}_{K}, \mathcal{H}_{K}, D_{K})$ 
is an even 
$QC^{\infty}$-spectral triple of 
spectral dimension 
$\dim_{S}(K)$. 
\end{theorem}

\begin{proof} 
By the definition of $D_{K}$ and 
Lemma \ref{lem:Ksummable}, 
$(\mathcal{A}_{K}, \mathcal{H}_{K}, D_{K})$ 
is an even spectral triple of spectral dimension 
$\dim_{S}(K)$. 
$(\mathcal{A}_{K}, \mathcal{H}_{K}, D_{K})$ is also of 
$QC^{\infty}$-class 
since we have $[|D_{K}| , T] = 0$ for 
an operator $T \in \mathcal{B}(\mathcal{H}_{K})$ of 
the direct sum of operators on $n$-cubes $f_{\bs{s}}(\gamma_{n})$. 
\end{proof}  

We next prove a nonvanishing property of 
the $K^{0}$-class of the Fredholm module $(\mathcal{H}_{K}, F_{K})$.

\begin{theorem}
\label{thm:suffnonvanish} 
Denote by $X_{1} , \dots , X_{k}$ 
the all connected components of 
$\displaystyle V \cup \bigcup_{s \in S} f_{s}(\gamma_{n})$. 
Then, if there is $X_{i}$ such that 
\[
\sharp (V_{0} \cap X_{i})
\neq 
\sharp (V_{1} \cap X_{i}), 
\]
the Connes-Chern character 
$\mathrm{Ch}_{\ast}(\mathcal{H}_{K}, F_{K}) 
	\in H_{\lambda}^{\mathrm{even}}(\mathcal{A}_{K})$
induces a non-zero additive map 
$K_{0}(C(V_{K})) \cong  K_{0}(\mathcal{A}_{K}) \to \mathbb{C}$ 
by the Connes pairing. 
In particular, 
$[\mathcal{H}_{K}, F_{K}]
\in K^{0}(C(V_{K}))$
is not trivial. 
\end{theorem}

\begin{proof}  
Set 
\[ 
d_{0} = 
\sharp (V_{0} \cap X_{i}), \quad 
d_{1} = \sharp (V_{1} \cap X_{i}) 
\]
and 
\begin{equation*}
p(x) = 
\begin{cases}
1 & x \in X_{i} \\ 
0 & \text{otherwise}
\end{cases} 
\end{equation*} 
for $x \in V_{K}$. 
Then,
$p$ is a continuous function 
and we have 
\begin{align*}
\ind (pF_{K}^{+}p : p\mathcal{H}_{K}^{+} \to p\mathcal{H}_{K}^{-}) 
&= 
\ind (pU_{n}p : p\ell^{2}(V_{0}) \to p\ell^{2}(V_{1}))  \\ 
&= 
d_{0} - d_{1} \neq 0. 
\end{align*} 
Therefore, we have 
$\mathrm{Ch}_{\ast}(\mathcal{H}_{K}, F_{K}) \neq 0$ 
on $K_{0}(C(V_{K}))$.

\end{proof} 

\begin{remark}
The assumption in 
Theorem \ref{thm:suffnonvanish} does not hold for some examples such as the Sierpinski carpet (see subsection \ref{subsec:SC}) 
and the $n$-cube $\gamma_{n}$. 
In these cases, the Connes-Chern character 
induces the $0$-map on $K_{0}(\mathcal{A}_{K})$. 
\end{remark}

\begin{remark}
\label{rmk:unitaryK}
As remarked in Remark \ref{rmk:unitary}, 
we can define a Fredholm module on $C(V)$ by using any unitary matrix $U$ 
instead of $U_{n}$. 
All properties in subsection \ref{subsec:FMST} hold 
without changing proofs in such a situation. 
\end{remark}

\subsection{Quantized differential form on self-similar sets}

Note that all similitudes on $\gamma_{n}$ form 
$f_{s}(\bs{x}) = r_{s}T_{s}\bs{x} + \bs{b}_{s}$ 
for an orthogonal matrix $T_{s} \in O(n)$ and $\bs{b}_{s} \in \mathbb{R}^{n}$. 
It is easy to calculate the quantum differential form 
$[F_{K} , x^{\alpha}]$ in the case for 
$\gamma_{n} = [0,1]^{n}$ and $T_{s} = E_{n}$ (for any $s \in S$), 
which is the direct sum of the matrix $d_{n}x^{\alpha}$; see Proposition 
\ref{prp:e_a}. 
We can also express  
$[F_{K}, x^{\alpha}]$ explicitly 
for the general case and 
show that they satisfy ``a variation'' of the Clifford relation.  

\begin{proposition}
\label{prp:qK}
We have 
\[
[F_{K} , x^{\alpha}][F_{K} , x^{\beta}] 
= 
\begin{cases}
-[F_{K} , x^{\beta}][F_{K} , x^{\alpha}] & \quad \alpha \neq \beta \\ 
\displaystyle -\bigoplus_{\bs{s} \in S^{\infty}}\dfrac{e_{\bs{s}}^{2}}{n}E_{2^{n}} & \quad  \alpha = \beta 
\end{cases}. 
\]
\end{proposition}

\begin{proof} 
Take an orthogonal matrix $T_{\bs{s}} = [t_{ij}]_{i,j} \in O(n)$ 
and a vector $\bs{b}_{\bs{s}} \in \mathbb{R}^{n}$ such that 
the image of 
the affine transformation 
$g_{\bs{s}}(\bs{x}) = e_{\bs{s}}T_{\bs{s}}\bs{x} + \bs{b}_{\bs{s}}$ 
of $[0,1]^{n}$ equals $f_{\bs{s}}(\gamma_{n})$  
and 
$g_{\bs{s}}(\bs{x})$  
preserves the numbering the vertices 
of $[0,1]^{n}$ and $f_{\bs{s}}(\gamma_{n})$.   
If we assume $\gamma_{n} = [0,1]^{n}$, we have $f_{\bs{s}} = g_{\bs{s}}$. 
Note that we have 
\begin{align*}
[F_{K} , x^{\alpha} ]|_{\mathcal{H}_{\bs{s}}} 
= 
\frac{1}{\sqrt{n}} 
\begin{bmatrix} 
 & -{}^{t}(\Delta_{n}x^{\alpha} \circ G_{n}) \\ 
\Delta_{n}x^{\alpha} \circ G_{n} & 
\end{bmatrix}. 
\end{align*}
Recall that  
$v_{2j} - v_{2i-1} 
= \pm e_{\bs{s}}T_{\bs{s}}\bs{e}_{k}$  
when $g_{\bs{s}}^{-1}(v_{2j})$ is connecting $g_{\bs{s}}^{-1}(v_{2i-1})$ by an edge 
of the $n$-cube $[0,1]^{n}$
parallel with $x^{k}$-direction  
and $T_{\bs{s}}\bs{e}_{k} = \sum_{\alpha = 1}^{n} t_{\alpha k}\bs{e}_{\alpha}$, 
and we have 
\[
[F_{K} , x^{\alpha} ]|_{\mathcal{H}_{\bs{s}}} 
= 
\frac{e_{\bs{s}}}{\sqrt{n}}
\sum_{j=1}^{n} t_{\alpha j}e_{(n)}^{j}. 
\]
Thus,
\begin{align*}
[F_{K} , x^{\alpha} ][F_{K} , x^{\beta} ]|_{\mathcal{H}_{\bs{s}}} 
&= 
\frac{e_{\bs{s}}^{2}}{n}
\left( \sum_{j=1}^{n} t_{\alpha j}e_{(n)}^{j} \right)
\left( \sum_{j=1}^{n} t_{\beta k}e_{(n)}^{k} \right) 
= 
\frac{e_{\bs{s}}^{2}}{n}
\sum_{j,k} t_{\alpha j}t_{\beta k}e_{(n)}^{j}e_{(n)}^{k} \\
&= 
\frac{e_{\bs{s}}^{2}}{n}
\sum_{j \neq k} t_{\alpha j}t_{\beta k}e_{(n)}^{j} e_{(n)}^{k}
- \frac{e_{\bs{s}}^{2}}{n} \sum_{j=1}^{n} t_{\alpha j}t_{\beta j} \\ 
&= 
\begin{cases}
\displaystyle 
\dfrac{e_{\bs{s}}^{2}}{n}
\sum_{j \neq k} t_{\alpha j}t_{\beta k}e_{(n)}^{j} e_{(n)}^{k} & (\alpha \neq \beta) \\ 
- \dfrac{e_{\bs{s}}^{2}}{n}E_{2^{n}} & (\alpha = \beta)
\end{cases}. 
\end{align*}
Therefore, we have
\[
[F_{K} , x^{\alpha}][F_{K} , x^{\beta}] 
= 
\begin{cases}
-[F_{K} , x^{\beta}][F_{K} , x^{\alpha}] & (\alpha \neq \beta) \\ 
\displaystyle -\bigoplus_{\bs{s} \in S^{\infty}}\dfrac{e_{\bs{s}}^{2}}{n}E_{2^{n}} 
	& (\alpha = \beta) 
\end{cases}. 
\]
 
\end{proof}

By Proposition \ref{prp:qK}, 
we get an explicit formula for 
$|[F_{K} , x^{1}] \cdots [F_{K} , x^{n}]|$. 

\begin{proposition}
\label{prp:absqf}
We have 
\[
|[F_{K} , x^{1}] \cdots [F_{K} , x^{n}]| = \bigoplus_{\bs{s} \in S^{\infty}} \frac{e_{\bs{s}}^{n}}{n^{n/2}}E_{2^{n}}.  
\]
\end{proposition}

\begin{proof}
Similar to the proof of Proposition \ref{prp:nDixmier}. 
\end{proof}

\begin{remark}
Set $\displaystyle e_{K}^{\alpha} = \bigoplus_{\bs{s} \in S^{\infty}} e_{(n)}^{\alpha}$, and we have the following Clifford relation 
\[
e_{K}^{\alpha}e_{K}^{\beta} 
= 
\begin{cases}
-e_{K}^{\beta}e_{K}^{\alpha} & (\alpha \neq \beta) \\ 
-\mathrm{id}_{H_{K}} & (\alpha \neq \beta)
\end{cases}. 
\]
Thus we can regard $e_{K}^{\alpha}$ as a $0$-$Q$-form in the sense of 
\cite{MR232790}. 
\end{remark}

\section{Dixmier traces} 
\label{sec:Dixmier}

In this section, we calculate the Dixmier trace of 
two operators. 
In general, the value for the second operator changes if  
the Fredholm operator $F_{n}$ changes 
to a different Fredholm operator. 

\subsection{Dixmier trace of \texorpdfstring{$|D_{K}|^{-p}$}{|DK|}}
\label{subsec:Dixv}

In this subsection, 
we calculate the Dixmier trace of $|D_{K}|^{-p}$. This is given by the residue at the pole of the zeta function 
$\zeta_{D_{K}}(s) = \mathrm{Tr}(|D_{K}|^{-s})$.

\begin{theorem}
\label{thm:Kzeta} 
For any $p \geq \dim_{S}(K)$, we have 
$|D_{K}|^{-p} \in \mathcal{L}^{(1,\infty)}(\mathcal{H}_{K})$ and 
\begin{equation*} 
\mathrm{Tr}_{\omega}(|D_{K}|^{-p}) 
= 
\begin{cases}
\displaystyle -2^{n}\left( \dim_{S}(K) \sum_{s=1}^{N} r_{s}^{\dim_{S}(K)} \log r_{s} \right)^{-1} & (p = \dim_{S}(K)) \\ 
0 & (p > \dim_{S}(K))
\end{cases}. 
\end{equation*}
Thus we have 
\begin{align*}
\mathrm{Tr}_{\omega}(f|D_{K}|^{-\dim_{S}(K)}) 
= -2^{n}\left( \dim_{S}(K) \sum_{s=1}^{N} r_{s}^{\dim_{S}(K)} \log r_{s} \right)^{-1}
\int_{K}f|_{K} \; d\Lambda 
\end{align*}
for any $f \in C(V_{K})$ by the Riesz-Markov-Kakutani representation theorem. 
Here, 
$\Lambda$ is 
the $\dim_{S}(K)$-dimensional Hausdorff probability measure of $K$. 

In particular, 
if all similarity ratios $r_{s}$ are equal, we have 
\[ 
\mathrm{Tr}_{\omega}(|D_{K}|^{-\dim_{S}(K)}) = \dfrac{2^{n}}{\log N}. 
\]
\end{theorem}

\begin{proof} 
By the proof of 
Lemma \ref{lem:Ksummable}, we have 
\begin{align*}
\mathrm{Tr}(|D_{K}|^{-p}) 
= 
2^{n} \sum_{j=0}^{\infty} \left( \sum_{s=1}^{N}r_{s}^{p}\right)^{j} 
= 
2^{n} \left( 1-\sum_{s=1}^{N}r_{s}^{p}\right)^{-1}.  
\end{align*}
Thus we have 
\begin{align*}
(z-1)\mathrm{Tr}(|D_{K}|^{-zp}) 
&= 
2^{n} \frac{z-1}{1-\sum_{s=1}^{N}r_{s}^{zp}} 
= 
2^{n} \frac{z-1}{\sum_{s=1}^{N} \left(r_{s}^{\dim_{S}(K)} - r_{s}^{zp}\right)} \\ 
&= 
2^{n} \left( \sum_{s=1}^{N} \dfrac{r_{s}^{\dim_{S}(K)} - r_{s}^{zp}}{z-1} \right)^{-1}. 
\end{align*}
and the following value 
\begin{align*}
\mathrm{Tr}_{\omega}(|D_{K}|^{-p})  
&= 
\lim_{z \to +1}(z-1)\mathrm{Tr}(|D_{K}|^{-zp}) \\  
&= 
2^{n} \left( \sum_{s=1}^{N} \lim_{z \to +1}\dfrac{r_{s}^{\dim_{S}(K)} - r_{s}^{zp}}{z-1} \right)^{-1}
\end{align*}
converges for  
$p \geq \dim_{S}(K)$. 
Finally, we get
\begin{align*}
\mathrm{Tr}_{\omega}(|D_{K}|^{-p})  
&= 
-2^{n}\left( \sum_{s=1}^{N}
\left. \dfrac{d}{dz}\right|_{z=1}r_{s}^{z\dim_{S}(K)} 
\right)^{-1} \\ 
&= 
-2^{n}\left( \dim_{S}(K) \sum_{s=1}^{N} r_{s}^{\dim_{S}(K)} \log r_{s} \right)^{-1}
\end{align*}
for 
$p = \dim_{S}(K)$ and 
\begin{align*}
\mathrm{Tr}_{\omega}(|D_{K}|^{-p})  = 0 
\end{align*}
for $p > \dim_{S}(K)$. 

\end{proof} 

\subsection{Dixmier trace of \texorpdfstring{$|[F_{K} , x^{1}] \cdots [F_{K} , x^{n}]|^{p}$}{quantized volume}}
\label{subsec:Dixqv}

In this subsection, 
we calculate the Dixmier trace of 
$|[F_{K} , x^{1}] \cdots [F_{K} , x^{n}]|^{p}$  
by using Proposition \ref{prp:absqf}. 

\begin{theorem}
\label{thm:Dix2}
We have 
$|[F_{K} , x^{1}][F_{K}, x^{2}] \cdots [F_{K}, x^{n}]|^{p}
\in \mathcal{L}^{(1,\infty)}(\mathcal{H}_{K})$ 
for any $p \geq \dfrac{1}{n}\dim_{S}(K)$. 
Moreover, we have 
\begin{align*}
&\; \mathrm{Tr}_{\omega}(|[F_{K} , x^{1}][F_{K}, x^{2}] \cdots [F_{K}, x^{n}]|^{p}) 
= 
\frac{1}{n^{np/2}}\mathrm{Tr}_{\omega}(|D_{K}|^{-np})  \\ 
=& 
\begin{cases}
\displaystyle 
\frac{-2^{n}}{n^{\dim_{S}(K)/2}}
\left( \dim_{S}(K) \sum_{s=1}^{N} r_{s}^{\dim_{S}(K)} \log r_{s} \right)^{-1}  
	& (p=\dfrac{1}{n}\dim_{S}(K)) \\ 
0 & (p>\dfrac{1}{n}\dim_{S}(K))  
\end{cases}. 
\end{align*}
Thus we have 
\begin{align*}
&\; \mathrm{Tr}_{\omega}(f|[F_{K} , x^{1}][F_{K}, x^{2}] \cdots [F_{K}, x^{n}]|^{\frac{1}{n}\dim_{S}(K)})  \\ 
=& 
\frac{-2^{n}}{n^{\dim_{S}(K)/2}}
\left( \dim_{S}(K) \sum_{s=1}^{N} r_{s}^{\dim_{S}(K)} \log r_{s} \right)^{-1} 
\int_{K}f|_{K} \; d\Lambda   \\ 
=& 
\frac{1}{n^{\dim_{S}(K)/2}}\mathrm{Tr}_{\omega}(|D_{K}|^{-\dim_{S}(K)}) 
\int_{K}f|_{K} \; d\Lambda   
\end{align*}
for any $f \in C(V_{K})$ by the Riesz-Markov-Kakutani representation theorem. 
Here, 
$\Lambda$ is 
the $\dim_{H}(K)$-dimensional Hausdorff probability measure of $K$. 
\end{theorem}

\begin{proof} 

By 
Proposition \ref{prp:absqf}, we have 
\begin{align*}
|[F_{K} , x^{1}][F_{K}, x^{2}] \cdots [F_{K}, x^{n}]|^{p}
&= \bigoplus_{\bs{s} \in S^{\infty}} 
	\frac{e_{\bs{s}}^{np}}{n^{np/2}}E_{2^{n}}. 
\end{align*}
Therefore, we get 
\begin{align*}
\mathrm{Tr}(|[F_{K} , x^{1}][F_{K}, x^{2}] 
\cdots [F_{K}, x^{n}]|^{p}) 
&= 
2^{n} \sum_{j=0}^{\infty} \sum_{(s_{1}, \dots s_{j}) \in S^{j}} 
\frac{1}{n^{np/2}}\prod_{k=1}^{j}r_{s_{k}}^{np} \\ 
&= 
\frac{2^{n}}{n^{np/2}} \sum_{j=0}^{\infty} \left( \sum_{s=1}^{N} r_{s}^{np} \right)^{j}
\end{align*}
and the following condition
\begin{align*}
\mathrm{Tr}(|[F_{K} , x^{1}][F_{K}, x^{2}] \cdots [F_{K}, x^{n}]|^{p}) < \infty 
\iff 
p > \frac{1}{n}\dim_{S}(K).
\end{align*}
If $p$ satisfies the above condition, the LHS can be written as 
\begin{align*}
\mathrm{Tr}(|[F_{K} , x^{1}][F_{K}, x^{2}] \cdots [F_{K}, x^{n}]|^{p})
&= 
\frac{2^{n}}{n^{np/2}} 
\left( 1- \sum_{s=1}^{N} r_{s}^{np} \right)^{-1}. 
\end{align*}
Therefore, the similar proof of Theorem \ref{thm:Kzeta} implies
\[ 
|[F_{K} , x^{1}][F_{K}, x^{2}] \cdots [F_{K}, x^{n}]|^{p} 
\in \mathcal{L}^{(1,\infty)}(\mathcal{H}_{K})
\] 
for 
$p \geq \dfrac{1}{n}\dim_{S}(K)$. 
Moreover, we get  
\begin{align*}
&\phantom{=} \mathrm{Tr}_{\omega}(|[F_{K} , x^{1}][F_{K}, x^{2}] \cdots [F_{K}, x^{n}]|^{p}) \\ 
&= 
\lim_{z \to +1} (z-1)
\mathrm{Tr}(|[F_{K} , x^{1}][F_{K}, x^{2}] \cdots [F_{K}, x^{n}]|^{zp}) \\ 
&= 
\frac{2^{n}}{n^{\dim_{S}(K)/2}}\left( \sum_{s=1}^{N} \lim_{z \to +1} \dfrac{r_{s}^{\dim_{S}(K)} - r_{s}^{z\dim_{S}(K)}}{z-1} \right)^{-1} \\ 
&= 
-\frac{2^{n}}{n^{\dim_{S}(K)/2}}
\left( \dim_{S}(K) \sum_{s=1}^{N} r_{s}^{\dim_{S}(K)} \log r_{s} \right)^{-1} 
\end{align*}
for  $p=\dfrac{1}{n}\dim_{S}(K)$ 
and 
\begin{align*}
\mathrm{Tr}_{\omega}(|[F_{K} , x^{1}][F_{K}, x^{2}] \cdots [F_{K}, x^{n}]|^{p})
= 
0  
\end{align*}
for  $p > \dfrac{1}{n}\dim_{S}(K)$. 

\end{proof}

\section{Examples}
\label{sec:exm}

In this section, 
we apply arguments in Section 
\ref{sec:FonK} and \ref{sec:Dixmier} to some examples. 

\subsection{Cantor dust} 
\label{subsec:CD}

The Cantor dust is a generalization of the middle third Cantor set 
to a higher dimension. 
Let $\mathcal{CD}_{n}$ be the Cantor dust defined on 
$\gamma_{n} = [0,1]^{n}$  
and the similitudes be
\[
f_{s}(\bs{x}) = \frac{1}{3}\bs{x} + \dfrac{2}{3}\sum_{\alpha = 1}^{n}a_{\alpha}\bs{e}_{\alpha} \quad 
	(\bs{x} \in \gamma_{n}, s = 0,1,2, \dots , 2^{n}-1).  
\]
Here, we write $a_{n}a_{n-1} \cdots a_{2}a_{1}$ as a number $s$ in binary and 
$\bs{e}_{\alpha}$ is the standard basis of $\mathbb{R}^{n}$.  
Since $\mathcal{CD}_{n}$ satisfies the open set condition, 
we have $\dim_{H}(\mathcal{CD}_{n}) = \dim_{S}(\mathcal{CD}_{n}) =  n\log_{3}2$. 
We also have $V_{\mathcal{CD}_{n}} = \mathcal{CD}_{n}$ since $\displaystyle V \subset \bigcup_{s=0}^{2^{n}-1}f_{s}(V)$. 
Then, we get 
\begin{align*}
\mathcal{A}_{\mathcal{CD}_{n}} = \mathrm{Lip}(\mathcal{CD}_{n}) 
\text{ and } 
C(V_{\mathcal{CD}_{n}}) = C(\mathcal{CD}_{n}). 
\end{align*}

Since all $f_{s}(\gamma_{n})$ are disconnected each other and  
also $\sharp (V_{0} \cap f_{1}(\gamma_{n})) = 1$ and 
$\sharp (V_{1} \cap f_{1}(\gamma_{n})) = 0$, 
the $K^{0}$-class of $(\mathcal{H}_{\mathcal{CD}_{n}}, F_{\mathcal{CD}_{n}})$ 
in $K^{0}(C(\mathcal{CD}_{n}))$
does not vanish by Theorem \ref{thm:suffnonvanish}.

\begin{theorem}
The Connes-Chern character 
\[ 
\mathrm{Ch}_{\ast}(\mathcal{H}_{\mathcal{CD}_{n}}, F_{\mathcal{CD}_{n}}) 
	\in H_{\lambda}^{\mathrm{even}}(\mathrm{Lip}(\mathcal{CD}_{n}))
\] 
induces a non-zero additive map 
$K_{0}(C(\mathcal{CD}_{n})) \to \mathbb{C}$. 
In particular, 
$[\mathcal{H}_{\mathcal{CD}_{n}}, F_{\mathcal{CD}_{n}}]$
is not trivial in 
$K^{0}(C(\mathcal{CD}_{n}))$. 
\end{theorem}

Since $\dim_{S}(\mathcal{CD}_{n}) = n\log_{3}2$, 
we also get the following results. 

\begin{corollary}
\begin{enumerate}
\item 
$(\mathcal{H}_{\mathcal{CD}_{n}}, F_{\mathcal{CD}_{n}})$ 
is a 
$([n\log_{3}2] + 1)$-summable even Fredholm module 
over $\mathrm{Lip}(\mathcal{CD}_{n})$. 
\item 
$(\mathrm{Lip}(\mathcal{CD}_{n}), \mathcal{H}_{\mathcal{CD}_{n}}, D_{\mathcal{CD}_{n}})$ 
is a 
$QC^{\infty}$-spectral triple 
of spectral dimension 
$n\log_{3}2$. 
\end{enumerate}
\end{corollary}

\begin{corollary}
We have the following. 
\begin{enumerate}
\item 
$\displaystyle 
\mathrm{Tr}(|D_{\mathcal{CD}_{n}}|^{-p}) 
= 
\frac{2^{n} \cdot 3^{p}}{3^{p} - 2^{n}}
$
for any $p > n\log_{3}2$. 
\item 
$\displaystyle 
\mathrm{Tr}_{\omega}(|D_{\mathcal{CD}_{n}}|^{-n\log_{3}2}) 
= 
\frac{2^{n}}{n\log 2}
$. 
\item 
$\displaystyle 
\mathrm{Tr}_{\omega}(f|D_{\mathcal{CD}_{n}}|^{-n\log_{3}2}) 
= 
\frac{2^{n}}{n\log 2}\int_{\mathcal{CD}_{n}}f \; d\Lambda
$ 
for any $f \in C(\mathcal{CD}_{n})$. 
Here, $\Lambda$ is the $(n\log_{3}2)$-dimensional Hausdorff 
probability measure of $\mathcal{CD}_{n}$. 
\end{enumerate}
\end{corollary}

\begin{corollary}
An operator 
$|[F_{\mathcal{CD}_{n}}, x^{1}][F_{\mathcal{CD}_{n}} , x^{2}] 
	\cdots [F_{\mathcal{CD}_{n}}, x^{n}]|^{\log_{3}2}$ 
is of $\mathcal{L}^{(1,\infty)}$-class and we have 
\begin{equation*}
\mathrm{Tr}_{\omega}(|[F_{\mathcal{CD}_{n}}, x^{1}][F_{\mathcal{CD}_{n}} , x^{2}] 
	\cdots [F_{\mathcal{CD}_{n}}, x^{n}]|^{\log_{3}2}) 
= 
\dfrac{2^{n}}{n^{(2+n\log_{3}2)/2}\log 2}.  
\end{equation*}
Thus we have 
\begin{align*}
\mathrm{Tr}_{\omega}(f|[F_{\mathcal{CD}_{n}}, x^{1}][F_{\mathcal{CD}_{n}} , x^{2}] 
	\cdots [F_{\mathcal{CD}_{n}}, x^{n}]|^{\log_{3}2}) 
&= 
\dfrac{2^{n}}{n^{(2+n\log_{3}2)/2}\log 2}
\int_{\mathcal{CD}_{n}}f \; d\Lambda  
\end{align*}
for any $f \in C(\mathcal{CD}_{n})$. 
Here, $\Lambda$ is the $(n\log_{3}2)$-dimensional Hausdorff 
probability measure of $\mathcal{CD}_{n}$. 
\end{corollary}

\subsection{Middle third Cantor set, revisited} 
\label{subsec:middle}

In this subsection, 
we focus on the middle third Cantor set $\mathcal{CS} = \mathcal{CD}_{1}$. 

First, we see a relationship between our Fredholm module and 
Connes' Fredholm module defined in \cite[Chapter IV. 3. $\varepsilon$]{MR1303779}. 
We recall Connes' Fredholm module $(H,F)$ on $C(\mathcal{CS})$. 
Let $I_{i,j} = (a_{i,j} , b_{i,j})$ ($i \in \mathbb{N}$, $j=1 , 2 , \dots , 2^i$)  
be open intervals in $[0,1]$ 
which are defined as 
\[
I_{1,1} = \left( \frac{1}{3}, \frac{2}{3}\right) \text{ and } 
I_{i+1,j} = \left( \frac{2b_{i,j-1} + a_{i,j}}{3} , \frac{b_{i,j-1} + 2 a_{i,j}}{3}\right). 
\]
Here, we set $b_{i,0} = 0$ and $a_{i,i+1} = 1$. 
The middle third Cantor set satisfies 
$\displaystyle 
\mathcal{CS} = [0,1] \setminus \bigcup_{i,j} I_{i,j}$. 
Connes defined 
\[
H = \bigoplus_{i,j} \ell^{2}(\{ a_{i,j} , b_{i,j} \}) \text{  and  } 
F = \bigoplus_{i,j} F_{1}. 
\]
Note that $H \oplus \ell^{2}(\{ 0,1 \}) \cong \mathcal{H}_{\mathcal{CS}}$ 
as Hilbert spaces. 

\begin{lemma}
\label{lem:composite}
Let $a < b < c$ be real numbers. 
We assume 
\[ 
[\ell^{2}(\{ a,b \}) , F_{1}] , [\ell^{2}(\{ b,c \}),F_{1}], [\ell^{2}(\{ a,c \}) , F_{1}]  
	\in K^{0}(C(\{ a,b,c \})) 
\] 
under homomorphisms $K^{0}(C(\{ a,b \})) \to K^{0}(C(\{ a,b,c \}))$, 
$K^{0}(C(\{ b,c \})) \to K^{0}(C(\{ a,b,c \}))$ and 
$K^{0}(C(\{ a,c \})) \to K^{0}(C(\{ a,b,c \}))$ 
defined by inclusions $\{ a,b \} \to \{ a,b,c \}$, 
$\{ b,c \} \to \{ a,b,c \}$ and 
$\{ a,c \} \to \{ a,b,c \}$, respectively. 
Then we have 
\[
[\ell^{2}(\{ a,b \}) , F_{1}] + [\ell^{2}(\{ b,c \}),F_{1}] 
= 
[\ell^{2}(\{ a,c \}) , F_{1}] \text{ in } K^{0}(C(\{ a,b,c \})). 
\]
\end{lemma}

\begin{proof}
Set $b=b_{1}=b_{2}$, 
$\{ a,b \} = \{ a,b_{1} \}$ and 
$\{ b,c \} = \{ b_{2},c \}$. 
We have 
\begin{align*}
[\ell^{2}(\{ a,b_{1} \}) , F_{1}] + [\ell^{2}(\{ b_{2},c \}),F_{1}]  
&= 
[\ell^{2}(\{ a,b_{1} \}) \oplus \ell^{2}(\{ b_{2},c \}) , F_{1} \oplus F_{1}] \\ 
&= 
\left[ \ell^{2}(\{ a,c \}) \oplus \ell^{2}(\{ b_{1},b_{2} \}) , 
\begin{bmatrix} & E_{2} \\ E_{2} & \end{bmatrix} \right]. 
\end{align*}
Here the $\mathbb{Z}_{2}$-grading operator of the last Fredholm module is 
defined by $\tilde{\epsilon} = \epsilon \oplus (-\epsilon)$. 
Set 
\[
T_{t} = \begin{bmatrix} F_{1}\cos t & \sin t \\ \sin t & -F_{1}\cos t \end{bmatrix} 
\]
on $\ell^{2}(\{ a,c \}) \oplus \ell^{2}(\{ b_{1},b_{2} \})$. 
Then we have 
$T_{t}\tilde{\epsilon} + \tilde{\epsilon}T_{t} = 0$, 
$T_{0} = F_{1} \oplus (-F_{1})$ and 
$T_{\pi/2} = \begin{bmatrix} & E_{2} \\ E_{2} & \end{bmatrix}$. 
Thus we have 
\begin{align*}
[\ell^{2}(\{ a,b_{1} \}) , F_{1}] + [\ell^{2}(\{ b_{2},c \}),F_{1}]  
&= 
\left[ \ell^{2}(\{ a,c \}) \oplus \ell^{2}(\{ b_{1},b_{2} \}) , 
F_{1} \oplus (-F_{1}) \right] \\ 
&= 
[ \ell^{2}(\{ a,c \}) , F_{1}] - [\ell^{2}(\{ b_{1},b_{2} \}) , F_{1}] \\ 
&=  
[ \ell^{2}(\{ a,c \}) , F_{1}]. 
\end{align*}
Here, the last equality is given by $b=b_{1}=b_{2}$. 
\end{proof}

By Lemma \ref{lem:composite}, we have 
\[
[H,F] + [\mathcal{H}_{\mathcal{CS}}, F_{\mathcal{CS}}] 
= [H_{\mathcal{CS}}, F_{\mathcal{CS}}]  + [\ell^{2}(\{ 0,1 \}) , F_{1}]. 
\]
Therefore we have $[H,F] = [\ell^{2}(\{ 0,1 \}) , F_{1}]$ 
in $K^{0}(C(\mathcal{CS}))$. 
On the other hand, if set 
\[
p_{k}(x) = 
\begin{cases}
1 & x \in [0,1/3^{k}] \cap \mathcal{CS} \\ 
0 & \text{otherwise}
\end{cases}
\]
for $x \in \mathcal{CS}$, then we get 
$\langle [H_{\mathcal{CS}}, F_{\mathcal{CS}}] , [p_{k}] \rangle 
= k$ 
and 
$\langle [\ell^{2}(\{ 0,1 \}) , F_{1}] , [p_{k}] \rangle 
= 1$
by the index pairing between $K$-homology and $K$-theory.  
Thus a pair 
$( [H_{\mathcal{CS}}, F_{\mathcal{CS}}] , [H,F] )$ 
is linearly independent on $\mathbb{Z}$ 
in $K^{0}(C(\mathcal{CS}))$. 

Second, we set similitudes 
\[
f_{1}(\bs{x}) = \frac{1}{3}\bs{x}, \quad  
f_{2}(\bs{x}) = \frac{1}{3}\bs{x} + \frac{2}{3}\bs{e}_{1}  
\]
for $\bs{x} \in \gamma_{2}$ 
and denote by $K$ the self-similar set defined by the IFS ($\gamma_{2} , \{ f_{1},f_{2} \}$). 
Then we get $K = \mathcal{CS} \times \{ 0 \}$ as sets. 
So the Fredholm module $(\mathcal{H}_{K} , F_{K})$ 
is a novel Fredholm module of the middle third Cantor set. 
Note that we have 
$V_{K} \neq K$ and 
$\left( \bigcup_{\bs{s} \in S^{\infty}}V_{\bs{s}}\right) \cap K \neq \emptyset$ 
in this case.

\begin{figure}[h]
  \includegraphics[width = 120mm]{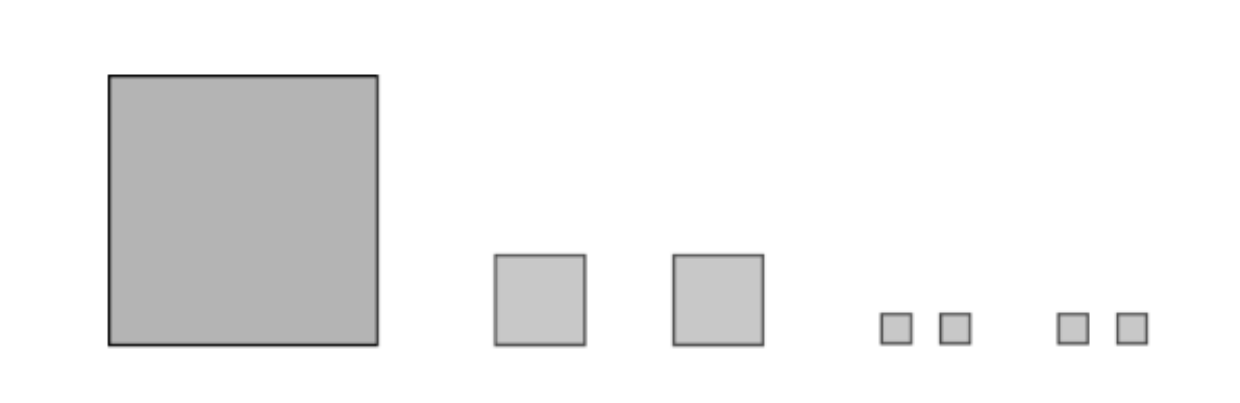}
\caption{The first $3$ steps of   construction of $K$. }
\end{figure}

\subsection{Sierpinski carpet and its higher dimensional analogue}
\label{subsec:SC}

The Sierpinski carpet is another generalization of 
the middle third Cantor set to a ``$2$-dimensional space''. 
The Menger sponge is also an analogue of 
the Sierpinski carpet but in a ``$3$-dimensional space''. 
In this subsection, we delve into such self-similar sets in $n$-dimensional space ($n \geq 2$). 
Let $S_{n} \subset \mathbb{N} \cup \{ 0 \}$ be the index set 
defined by 
\begin{align*} 
S_{n} = \{ s \in \mathbb{N} \cup \{ 0 \} \,;\; &0 \leq s \leq 3^{n}-1 
\text{ and at most one} \\ &\text{ of digits equals } 1   
\text{ in ternary expression of }  s 
  \}. 
\end{align*} 
For example, for $n = 2, 3$, we have 
$S_{2} = \{ 0,1,2,3,5,6,7,8 \}$ and 
\[ 
S_{3} = S_{2} \cup \{ 9,11,15,17,18,19,20,21,23,24,25,26 \}. 
\]
Define similitudes $f_{s} : \gamma_{n} \to \gamma_{n}$ for $s \in S_{n}$ 
by 
\[
f_{s}(\bs{x}) = \frac{1}{3}\bs{x} + \frac{1}{3}\sum_{\alpha = 1}^{n}a_{\alpha}\bs{e}_{\alpha}. 
\]
Here, we use a number $s$ to express $a_{n}a_{n-1} \cdots a_{2}a_{1}$ in ternary. 
We write $\mathcal{SC}_{n}$ as the self-similar set on 
the IFS $(\gamma_{n}, S_{n}, \{ f_{s} \}_{s \in S_{n}})$. 
When $n = 2$ and $3$, $\mathcal{SC}_{2}$ is the Sierpinski carpet and 
$\mathcal{SC}_{3}$ is the Menger sponge. 
Since $\mathcal{CD}_{n}$ satisfies the open set condition, 
we have $\dim_{H}(\mathcal{CD}_{n}) = \dim_{S}(\mathcal{CD}_{n}) = \log_{3}(\sharp S_{n})
= \log_{3} (2^{n-1}(n+2))$.  
We have $V_{\mathcal{SC}_{n}} = \mathcal{SC}_{n}$ since 
we have $\displaystyle V \subset \bigcup_{s \in S_{n}}f_{s}(V)$. 
Then, we get 
\begin{align*}
\mathcal{A}_{\mathcal{SC}_{n}} = \mathrm{Lip}(\mathcal{SC}_{n}), \quad 
C(V_{\mathcal{SC}_{n}}) = C(\mathcal{SC}_{n}). 
\end{align*}

\begin{figure}[h]
  \includegraphics[width = 120mm]{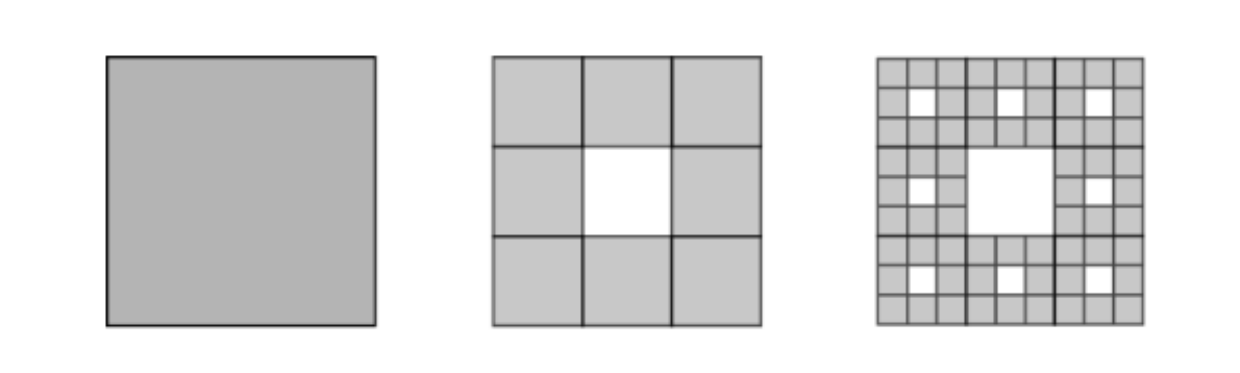}
\caption{The first $3$ steps of  construction of $\mathcal{SC}_{2}$. }
\end{figure}

Since $X = V \cup \bigcup_{s \in S_{n}} f_{s}(\gamma_{n})$ is 
connected,  we have 
$\sharp (V_{0} \cap X) = \sharp (V_{1} \cap X)$; the assumption in  
Theorem \ref{thm:suffnonvanish} does not hold. 

\begin{remark}
\label{rmk:SC2}
The Sierpinski carpet $\mathcal{SC}_{2}$ is a compact set in $\mathbb{R}^{2}$. 
We have $K_{0}(C(\mathcal{SC}_{2})) = \mathbb{Z}$ which is 
generated by (matrix valued) constant functions on $\mathcal{SC}_{2}$, and the index pairing between $K$-theory and $K$-homology 
induces the $0$-map $K_{0}(C(\mathcal{SC}_{2})) \to \mathbb{Z}$.  
Therefore we get $[\mathcal{H}_{\mathcal{SC}_{2}}, F_{\mathcal{SC}_{2}}] = 0$ 
in $K^{0}(C(\mathcal{SC}_{2}))$ by 
\cite[Theorem 7.5.5]{MR1817560}. 

On the other hand, 
we can construct a non-trivial Fredholm module corresponding to 
the Sierpinski carpet in a manner similar to the construction shown in subsection \ref{subsec:middle}. 
Define $z : \gamma_{1} \to \gamma_{1}$ by $z(t) = \dfrac{1}{3}t$ and $\tilde{f}_{s} = (f_{s}, z) : \gamma_{3} \to \gamma_{3}$ 
for $s \in S_{2}$. Then we get 
a new IFS $(\gamma_{3} , S_{2} , \{ \tilde{f}_{s} \}_{s \in S_{2}})$. 
Denote by $\widetilde{\mathcal{SC}}_{2}$ 
the self-similar set on the new IFS,  
and we get $\widetilde{\mathcal{SC}}_{2} = \mathcal{SC}_{2} \times \{ 0 \}$. 
The corresponding Fredholm module 
$(\mathcal{H}_{\widetilde{\mathcal{SC}}_{2}} , F_{\widetilde{\mathcal{SC}}_{2}})$ 
represents a non-trivial element in $K^{0}(C(V_{\widetilde{\mathcal{SC}}_{2}}))$. 
\end{remark}

\begin{remark}
The construction of IFS in Remark \ref{rmk:SC2} 
can be extended to 
general cases: let $(\gamma_{n},S,\{ f_{s} \}_{s \in S})$ be 
an IFS and $K$ its self-similar set. 
Then $(\gamma_{n+1}, S , \{ (f_{s} , z) \}_{s \in S})$ 
is a new IFS and the corresponding self-similar set 
denote by $\widetilde{K}$ satisfies $\widetilde{K} = K \times \{ 0 \}$ 
and $[\mathcal{H}_{\widetilde{K}} , F_{\widetilde{K}}] \neq 0$ 
in $K^{0}(C(V_{\widetilde{K}}))$. 
\end{remark}

Since $\dim_{S}(\mathcal{SC}_{n})  
= \log_{3} (2^{n-1}(n+2))$,  
we get the following results.

\begin{corollary}
\begin{enumerate}
\item 
$(\mathcal{H}_{\mathcal{SC}_{n}}, F_{\mathcal{SC}_{n}})$ 
is a 
$([\log_{3} (2^{n-1}(n+2))] + 1)$-summable even Fredholm module 
over $\mathrm{Lip}(\mathcal{SC}_{n})$. 
\item  
$(\mathrm{Lip}(\mathcal{SC}_{n}), \mathcal{H}_{\mathcal{SC}_{n}}, D_{\mathcal{SC}_{n}})$ 
is a 
$QC^{\infty}$-spectral triple 
of spectral dimension \\ 
$\log_{3} (2^{n-1}(n+2))$. 
\end{enumerate}
\end{corollary}

\begin{corollary}
We have the following. 
\begin{enumerate}
\item 
$\displaystyle 
\mathrm{Tr}(|D_{\mathcal{SC}_{n}}|^{-p}) 
= 
\frac{2^{n} \cdot 3^{p}}{3^{p} - 2^{n-1}(n+2)}
$
for any $p > \log_{3} (2^{n-1}(n+2))$. 
\item 
$\displaystyle 
\mathrm{Tr}_{\omega}(|D_{\mathcal{SC}_{n}}|^{-\log_{3}(2^{n-1}(n+2))}) 
= 
\frac{2^{n}}{\log (2^{n-1}(n+2))}
$. 
\item 
$\displaystyle 
\mathrm{Tr}_{\omega}(f|D_{\mathcal{SC}_{n}}|^{-\log_{3}(2^{n-1}(n+2))}) 
= 
\frac{2^{n}}{\log (2^{n-1}(n+2))}\int_{\mathcal{SC}_{n}}f \; d\Lambda
$ 
for any $f \in C(\mathcal{SC}_{n})$. 
Here, $\Lambda$ is the $(\log_{3}(2^{n-1}(n+2)))$-dimensional Hausdorff 
probability measure of $\mathcal{SC}_{n}$. 
\end{enumerate}
\end{corollary}

\begin{corollary}
For $d = \dfrac{1}{n}\log_{3}(2^{n-1}(n+2))$, 
we have 
\[ 
|[F_{\mathcal{SC}_{n}}, x^{1}][F_{\mathcal{SC}_{n}} , x^{2}] 
	\cdots [F_{\mathcal{SC}_{n}}, x^{n}]|^{d} \in \mathcal{L}^{(1,\infty)}(\mathcal{H}_{\mathcal{SC}_{n}})  
\] 
and  
\begin{equation*}
\mathrm{Tr}_{\omega}(|[F_{\mathcal{SC}_{n}}, x^{1}][F_{\mathcal{SC}_{n}} , x^{2}] 
	\cdots [F_{\mathcal{SC}_{n}}, x^{n}]|^{d}) 
= 
\dfrac{2^{n}}{n^{nd/2}\log (2^{n-1}(n+2))}  . 
\end{equation*}
Thus we have 
\begin{align*}
\mathrm{Tr}_{\omega}(f|[F_{\mathcal{CD}_{n}}, x^{1}][F_{\mathcal{CD}_{n}} , x^{2}] 
	\cdots [F_{\mathcal{CD}_{n}}, x^{n}]|^{d}) 
&= 
\dfrac{2^{n}}{n^{nd/2}\log (2^{n-1}(n+2))}
\int_{\mathcal{SC}_{n}}f \; d\Lambda  
\end{align*}
for any $f \in C(\mathcal{SC}_{n})$. 
Here, $\Lambda$ is the $(\log_{3}(2^{n-1}(n+2)))$-dimensional Hausdorff 
probability measure of $\mathcal{SC}_{n}$. 
\end{corollary}

\subsection{With rotations} 
\label{subsec:rotation}

Let $R = 
\begin{bmatrix} \cos \theta & -\sin \theta \\ \sin \theta & \cos \theta \end{bmatrix}$ the be a rotation matrix. 
Let also $f_{1},f_{2},f_{3},f_{4}$ be four similitudes defined by 
\begin{align*}
f_{s}(\bs{x}) 
= \frac{1}{2\sqrt{2}}R\left(\bs{x} - \frac{1}{2}\begin{bmatrix}1 \\ 1 \end{bmatrix}\right) 
	+ \bs{b}_{s}. 
\end{align*}
Here, we set 
\begin{align*}
\bs{b}_{1} = \frac{1}{4}\begin{bmatrix}1 \\ 1 \end{bmatrix}, \quad 
\bs{b}_{2} = \frac{1}{4}\begin{bmatrix}3 \\ 1 \end{bmatrix}, \quad 
\bs{b}_{3} = \frac{1}{4}\begin{bmatrix}1 \\ 3 \end{bmatrix}, \quad 
\bs{b}_{4} = \frac{1}{4}\begin{bmatrix}3 \\ 3 \end{bmatrix}. 
\end{align*}
The IFS $(\gamma_{2} , \{f_{1},f_{2},f_{3},f_{4}\})$ 
is defined by using 
a rotation of angle $\theta$. 
We get the self-similar set $K$ on the IFS $(\gamma_{2} , \{f_{1},f_{2},f_{3},f_{4}\})$ 
that satisfies the open set condition.  
Then we have 
$V_{K} \neq K$ and 
$\left( \bigcup_{\bs{s} \in \{ 1,2,3,4 \}^{\infty}}V_{\bs{s}} \right) \cap K = \emptyset$. 
Since $\{ (0,0) \}$ is a connected component of 
$\displaystyle V \cup \bigcup_{s \in \{1,2,3,4\}} f_{s}(\gamma_{2})$, 
the Fredholm module $(\mathcal{H}_{K}, F_{K})$ defines 
a non-trivial element in $K^{0}(C(V_{K}))$. 

\begin{figure}[h]
  \includegraphics[width = 120mm]{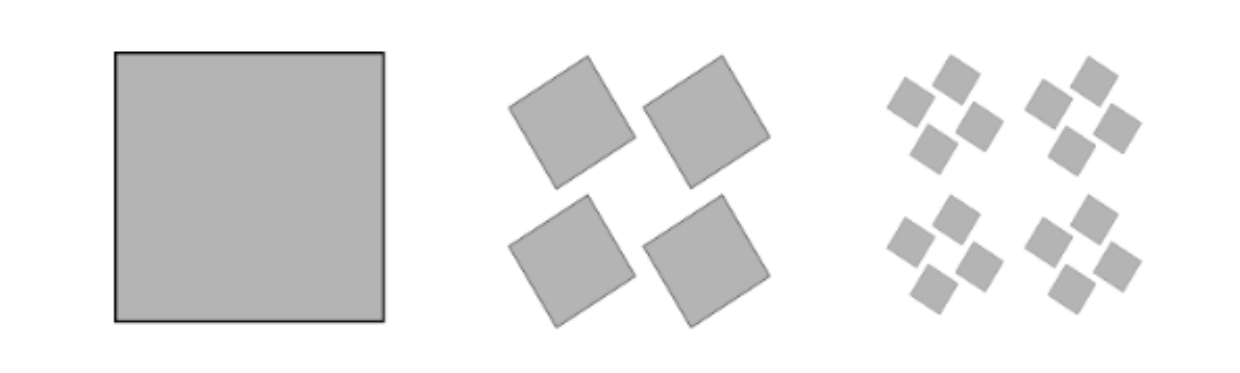}
\caption{The first $3$ steps of  construction of $K$. }
\end{figure}

Since $\dim_{S}(K) = \log_{2\sqrt{2}}4 = \dfrac{4}{3}$, 
we get the following results. 

\begin{corollary}
\begin{enumerate}
\item 
$(\mathcal{H}_{K}, F_{K})$ 
is a 
$2$-summable even Fredholm module 
over $\mathcal{A}_{K}$. 
\item 
$(\mathcal{A}_{K}, \mathcal{H}_{K}, D_{K})$ 
is a 
$QC^{\infty}$-spectral triple 
of spectral dimension 
$\dfrac{4}{3}$. 
\end{enumerate}
\end{corollary}

\begin{corollary}
We have the following. 
\begin{enumerate}
\item 
$\displaystyle 
\mathrm{Tr}(|D_{K}|^{-p}) 
= 
\frac{4}{2^{3p/2} - 4}
$
for any $p > \dfrac{4}{3}$. 
\item 
$\displaystyle 
\mathrm{Tr}_{\omega}(|D_{K}|^{-4/3}) 
= 
\frac{2}{\log 2}
$. 
\item 
$\displaystyle 
\mathrm{Tr}_{\omega}(f|D_{K}|^{-4/3}) 
= 
\frac{2}{\log 2}\int_{K}f|_{K} \; d\Lambda
$ 
for any $f \in C(V_{K})$. 
Here, $\Lambda$ is the $4/3$-dimensional Hausdorff 
probability measure of $K$. 
\end{enumerate}
\end{corollary}

By Proposition \ref{prp:df}, the quantized differential forms $[F_{K}, x^{\alpha}]$ ($\alpha = 1,2$)
are given as 
\begin{align*}
[F_{K}, x^{1}] 
&= 
\bigoplus_{j=0}^{\infty} \bigoplus_{\bs{s} \in S^{\times j}}
\frac{e_{\bs{s}}}{\sqrt{2}} 
\begin{bmatrix} 
0 & 0 & \cos j\theta & -\sin j\theta \\ 
0 & 0 & -\sin j\theta & -\cos j\theta \\ 
-\cos j\theta & \sin j\theta & 0 & 0 \\ 
\sin j\theta & \cos j\theta & 0 & 0 
\end{bmatrix}, \\
[F_{K},x^{2}] 
&= \bigoplus_{j=0}^{\infty} \bigoplus_{\bs{s} \in S^{\times j}}
\frac{e_{\bs{s}}}{\sqrt{2}} 
\begin{bmatrix} 
0 & 0 & \sin j\theta & \cos j\theta \\ 
0 & 0 & \cos j\theta & -\sin j\theta \\ 
-\sin j\theta & -\cos j\theta & 0 & 0 \\ 
-\cos j\theta & \sin j\theta & 0 & 0 
\end{bmatrix}. 
\end{align*}
Thus, we have 
\begin{align*}
|[F_{K}, x^{1}][F_{K}, x^{2}]| 
&= 
\bigoplus_{\bs{s} \in S^{\infty}} \frac{e_{\bs{s}}^{2}}{2}E_{4}. 
\end{align*}
This implies  
\begin{corollary}
An operator 
$|[F_{K}, x^{1}][F_{K} , x^{2}]|^{2/3}$ 
is of $\mathcal{L}^{(1,\infty)}$-class and we have 
\begin{equation*}
\mathrm{Tr}_{\omega}(|[F_{K}, x^{1}][F_{K} , x^{2}]|^{2/3}) 
= 
\dfrac{\sqrt[3]{2}}{\log 2}.   
\end{equation*}
Thus we have 
\begin{align*}
\mathrm{Tr}_{\omega}(f|[F_{K}, x^{1}][F_{K} , x^{2}]|^{2/3}) 
&= 
\dfrac{\sqrt[3]{2}}{\log 2} 
\int_{K}f|_{K} \; d\Lambda  
\end{align*}
for any $f \in C(V_{K})$. 
Here, $\Lambda$ is the $\dfrac{4}{3}$-dimensional Hausdorff 
probability measure of $K$. 
\end{corollary}

\subsection{Without the open set condition} 
\label{subsec:nonOSC}

In this subsection, 
we present an example of a self-similar set that does not satisfy 
the open set condition. 
In this case, we can detect the similarity dimension by using our Fredholm module 
but not detect the Hausdorff dimension explicitly. 

Let $(\gamma_{2} , S = \{ 1,2,3,4,5 \} , \{ f_{s} \}_{s \in S})$ 
be an IFS defined to be  
\begin{align*}
f_{1}(\bs{x}) &= \frac{1}{3}\bs{x}, \quad 
f_{2}(\bs{x}) = \frac{1}{3}\bs{x} + \frac{2}{3}\bs{e}_{1}, \quad 
f_{3}(\bs{x}) = \frac{1}{3}\bs{x} + \frac{2}{3}\bs{e}_{2}, \\  
f_{4}(\bs{x}) &= \frac{1}{3}\bs{x} + \frac{2}{3}\bs{e}_{1} + \frac{2}{3}\bs{e}_{2}, \quad 
f_{5}(\bs{x}) = \frac{2}{3}\bs{x} + \frac{1}{6}\bs{e}_{1} + \frac{1}{6}\bs{e}_{2}. 
\end{align*}
Note that the IFS does not satisfy the open set condition. 
Let $K$ be the self-similar set on the IFS. 
Since 
we have $\displaystyle V \subset \bigcup_{s =1}^{5}f_{s}(V)$, 
we have $V_{K} = K$. 
The similarity dimension $s = \dim_{S}(K)$ of $K$ is given by 
the following identity 
\[
4 \cdot \left( \frac{1}{3} \right)^{s} + \left( \frac{2}{3} \right)^{s} = 1.   
\]
We can easily check $1 < s < 2$.  

\begin{corollary}
\begin{enumerate}
\item 
$(\mathcal{H}_{K}, F_{K})$ 
is a 
$2$-summable even Fredholm module 
over $\mathrm{Lip}(K)$. 
\item 
$(\mathrm{Lip}(K), \mathcal{H}_{K}, D_{K})$ 
is a 
$QC^{\infty}$-spectral triple 
of spectral dimension 
$s$. 
\end{enumerate}
\end{corollary}

\begin{corollary}
We have the following. 
\begin{enumerate}
\item 
$\displaystyle 
\mathrm{Tr}(|D_{K}|^{-p}) 
= 
\frac{4 \cdot 3^{p}}{3^{p} - 2^{p} - 4}
$
for any $p > s$. 
\item 
$\displaystyle 
\mathrm{Tr}_{\omega}(|D_{K}|^{-\dim_{S}(K)}) 
= 
\frac{4 \cdot 3^{s}  }{  3^{s} s \log 3 -   2^{s} s \log 2} 
$. 
\item 
$\displaystyle 
\mathrm{Tr}_{\omega}(f|D_{K}|^{-\dim_{S}(K)}) 
= 
\frac{4 \cdot 3^{s}  }{ 3^{s} s \log 3 -  2^{s} s \log 2} 
\int_{K}f \; d\Lambda
$ 
for any $f \in C(K)$. 
Here, $\Lambda$ is the $\dim_{H}(K)$-dimensional Hausdorff 
probability measure of $K$. 
\end{enumerate}
\end{corollary}

\begin{corollary}
An operator 
$|[F_{K}, x^{1}][F_{K} , x^{2}]|^{s/2}$ 
is of $\mathcal{L}^{(1,\infty)}$-class and we have 
\begin{equation*}
\mathrm{Tr}_{\omega}(|[F_{K}, x^{1}][F_{K} , x^{2}]|^{d}) 
= 
\frac{2^{2-s/2} \cdot 3^{s} }{ 3^{s} s \log 3 -  2^{s} s \log 2}.   
\end{equation*}
Thus we have 
\begin{align*}
\mathrm{Tr}_{\omega}(f|[F_{K}, x^{1}][F_{K} , x^{2}]|^{d}) 
&= 
\frac{2^{2-s/2} \cdot 3^{s} }{ 3^{s} s \log 3 -  2^{s} s \log 2}    
\int_{K}f \; d\Lambda  
\end{align*}
for any $f \in C(K)$. 
Here, $\Lambda$ is the $\dim_{H}(K)$-dimensional Hausdorff 
probability measure of $K$. 
\end{corollary}



\begin{funding}
Seto was partially supported by JSPS KAKENHI Grant Number  21K13795.  
\end{funding}


\bibliographystyle{emss}
\bibliography{Fredholm_module_on_n-cube}









\end{document}